\documentclass[thmsa,12pt]{article}
\usepackage{amsmath}
\usepackage{amssymb}
\usepackage{sw20lart}

\input{tcilatex}

\begin{document}

\begin{center}
\medskip

\textbf{\textit{G}}$^{+}$ \textbf{Method and \v{C}ern\'{y} Conjecture}

\bigskip

\textbf{Udrea\ P\u{a}un}

\bigskip
\end{center}

\noindent We revisit and extend the $G^{+}$ method --- some results for the
sum-positive matrices are improved; the $a$-sum-positive matrices are
introduced, with theory and applications; some combinatorial problems are
proposed. The \v{C}ern\'{y} conjecture is a combinatorial conjecture for the
synchronizable deterministic finite automata on the length of reset words.
Besides this conjecture, we consider other two combinatorial conjectures,
the extended \v{C}ern\'{y} conjecture and more extended \v{C}ern\'{y}
conjecture --- the extended \v{C}ern\'{y} conjecture holds if and only if
the \v{C}ern\'{y} conjecture holds. To prove the \v{C}ern\'{y} conjecture
for classes of automata, several methods are used, two of them being the
compression method and extension method --- these two are submethods of the $%
G^{+}$ method, and, to illustrate this thing, we give two examples.

\medskip


\medskip


\begin{center}
\medskip

\textbf{1. \textit{G}}$^{+}$ \textbf{METHOD REVISITED AND\ EXTENDED}

\medskip
\end{center}

In this section, we revisit and extend the $G^{+}$ method [28], the second $%
G $-type method (for the study of products of nonnegative matrices) --- we
improve some results for the sum-positive matrices; we introduce the $a$%
-sum-positive matrices, with theory and two applications, one of them being
to finite Markov chains and the other one being in the reliability theory;
we propose some combinatorial problems.

\smallskip

Set 
\begin{equation*}
\left\langle m\right\rangle =\left\{ 1,2,...,m\right\} \text{ (}m\in \mathbb{%
N},\text{ }m\geq 1\text{)},
\end{equation*}
\begin{equation*}
\left\langle \left\langle m\right\rangle \right\rangle =\left\{
0,1,...,m\right\} \text{ (}m\in \mathbb{N}\text{)},
\end{equation*}

\begin{equation*}
N_{m,n}=\left\{ P\left| \text{ }P\text{ is a nonnegative }m\times n\text{
matrix}\right. \right\} \text{,}
\end{equation*}
\begin{equation*}
S_{m,n}=\left\{ P\left| \text{ }P\text{ is a stochastic }m\times n\text{
matrix}\right. \right\}
\end{equation*}

\noindent --- here, a stochastic matrix is a row stochastic matrix ---,

\begin{equation*}
N_{n}=N_{n,n}\text{,}
\end{equation*}

\begin{equation*}
S_{n}=S_{n,n}.
\end{equation*}

\smallskip

Let $Z$ be an $m\times n$ matrix. $Z_{ij}$ or, if confusion can arise, $%
Z_{i\rightarrow j}$ denotes the entry $\left( i,j\right) $ of matrix $Z,$ $%
\forall i\in \left\langle m\right\rangle ,$ $\forall j\in \left\langle
n\right\rangle .$

\smallskip

\textbf{Definition 1.1.} Let $P=\left( P_{ij}\right) \in N_{m,n}$ ($i\in
\left\langle m\right\rangle $ and $j\in \left\langle n\right\rangle $). $%
\left\langle m\right\rangle $ and $\left\langle n\right\rangle $ are called
the \textit{index sets of} $P;$ $\left\langle m\right\rangle $ is called the 
\textit{row index set of} $P$ while $\left\langle n\right\rangle $ is called
the \textit{column index set of} $P.$ If $P\in N_{n},$ $\left\langle
n\right\rangle $ is called the \textit{index set of} $P.$

\smallskip

Let $P=\left( P_{ij}\right) \in N_{m,n}$. Let $\emptyset \neq U\subseteq
\left\langle m\right\rangle $ and $\emptyset \neq V\subseteq \left\langle
n\right\rangle $. Set the matrices (these are submatrices of $P$) 
\begin{equation*}
P_{U}=\left( P_{ij}\right) _{i\in U,j\in \left\langle n\right\rangle },\text{
}P^{V}=\left( P_{ij}\right) _{i\in \left\langle m\right\rangle ,j\in V},%
\text{ and }P_{U}^{V}=\left( P_{ij}\right) _{i\in U,j\in V}.
\end{equation*}

\smallskip

\textbf{Definition 1.2.} (See, \textit{e.g.}, [35, p. 80].) Let $P\in
N_{m,n}.$ We say that $P$ is a \textit{row-allowable matrix} if it has at
least one positive entry in each row.

\smallskip

Below we define a central notion of this article.

\smallskip

\textbf{Definition 1.3} ([28]). Let $P\in N_{m,n}$. Let $\emptyset \neq
U\subseteq \left\langle m\right\rangle $ and $\emptyset \neq V\subseteq
\left\langle n\right\rangle $. We say that $P$ is \textit{sum-positive on }$%
U\times V$ if 
\begin{equation*}
\sum\limits_{j\in V}P_{ij}>0,\text{ }\forall i\in U
\end{equation*}
\noindent (equivalently, if $P_{U}^{V}$ is a row-allowable matrix).

\smallskip

For the sum-positive matrices on $U\times V,$ $\emptyset \neq U\subseteq
\left\langle m\right\rangle $, $\emptyset \neq V\subseteq \left\langle
n\right\rangle ,$ we will use the generic name ``sum-positive matrices''.

\smallskip

\textbf{Definition 1.4.} (See, \textit{e.g.}, [18, p. 57] and [35, p. 140].)
Let $P\in N_{m,n},$ $m,$ $n\geq 1$. We say that $P$ is a \textit{Markov
matrix} if it has at least one positive column (equivalently, if it is
sum-positive on $\left\langle m\right\rangle \times \left\{ j\right\} $ for
some $j\in \left\langle n\right\rangle $).

\smallskip

Let $P\in N_{m,n},$ $m,$ $n\geq 1$. Set

\begin{equation*}
\mu \left( P\right) =\max_{j\in \left\langle n\right\rangle }\left(
\min_{i\in \left\langle m\right\rangle }P_{ij}\right) .
\end{equation*}

\noindent $\mu \left( P\right) $ is called the \textit{Markov ergodicity
coefficient} of $P$ if $P\in S_{m,n}$ (see, \textit{e.g.}, [18, p. 56]).

\smallskip

\textbf{Remark 1.5.} Let $P\in N_{m,n}$ ($m,$ $n\geq 1$). $P$ is a Markov
matrix $\Longleftrightarrow $ $\mu \left( P\right) >0$ --- thus we have a
reason for the name ``Markov matrix''.

\smallskip

The Markov stochastic matrices are of interest, \textit{e.g.}, in the finite
Markov chain theory (see, \textit{e.g.}, [28, pp. 136--137]), in the
deterministic finite automaton theory (see Section 2...)...

\smallskip

Let $\emptyset \neq U\subseteq \left\langle m\right\rangle $ and $\emptyset
\neq V\subseteq \left\langle n\right\rangle $. Set (see [28]) 
\begin{equation*}
G_{U,V}^{+}=G_{U,V}^{+}\left( m,n\right) =\left\{ P\left| \text{ }P\in
S_{m,n}\text{ and is sum-positive on }U\times V\right. \right\}
\end{equation*}

\noindent and 
\begin{equation*}
\overline{G}_{U,V}^{+}=\overline{G}_{U,V}^{+}\left( m,n\right) =\left\{
P\left| \text{ }P\in N_{m,n}\text{ and is sum-positive on }U\times V\right.
\right\}
\end{equation*}

\noindent --- $G_{U,V}^{+}$ and $\overline{G}_{U,V}^{+}$ are two important
sets of sum-positive matrices; obviously, 
\begin{equation*}
G_{U,V}^{+}\subseteq S_{m,n},\text{ }\overline{G}_{U,V}^{+}\subseteq N_{m,n},%
\text{ and }G_{U,V}^{+}\subseteq \overline{G}_{U,V}^{+}.
\end{equation*}

\smallskip

Theorem 2.6 from [28] is a basic result for the sum-positive matrices, a
if-then theorem/result --- below we give a better theorem, a if-and-only-if
theorem.

\smallskip

\textbf{THEOREM 1.6.} (i) \textit{Let} $P_{1}\in N_{m_{1},m_{2}}$ \textit{%
and }$P_{2}\in N_{m_{2},m_{3}}.$\textit{\ Let }$\emptyset \neq
U_{1}\subseteq \left\langle m_{1}\right\rangle $ \textit{and} $\emptyset
\neq U_{3}\subseteq \left\langle m_{3}\right\rangle .$ \textit{Then} 
\begin{equation*}
P_{1}P_{2}\in \overline{G}_{U_{1},U_{3}}^{+}\Longleftrightarrow \exists
U_{2},\text{ }\emptyset \neq U_{2}\subseteq \left\langle m_{2}\right\rangle ,
\end{equation*}
\noindent \textit{such that} 
\begin{equation*}
P_{1}\in \overline{G}_{U_{1},U_{2}}^{+}\text{\hspace{0.2cm}\textit{and} }%
P_{2}\in \overline{G}_{U_{2},U_{3}}^{+}.
\end{equation*}

\smallskip

(ii) (\textit{a generalization of }(i)) \textit{Let }$P_{1}\in
N_{m_{1},m_{2}},$\textit{\ }$P_{2}\in N_{m_{2},m_{3}},$ $...,$ $P_{n}\in
N_{m_{n},m_{n+1}}.$\textit{\ Let }$\emptyset \neq U_{1}\subseteq
\left\langle m_{1}\right\rangle $ \textit{and} $\emptyset \neq
U_{n+1}\subseteq \left\langle m_{n+1}\right\rangle .$ \textit{Then} 
\begin{equation*}
P_{1}P_{2}...P_{n}\in \overline{G}_{U_{1},U_{n+1}}^{+}\Longleftrightarrow
\exists U_{2},U_{3},...,U_{n},
\end{equation*}
\begin{equation*}
\emptyset \neq U_{2}\subseteq \left\langle m_{2}\right\rangle ,\text{ }%
\emptyset \neq U_{3}\subseteq \left\langle m_{3}\right\rangle ,\text{ }...,%
\text{ }\emptyset \neq U_{n}\subseteq \left\langle m_{n}\right\rangle ,
\end{equation*}
\noindent \textit{such that} 
\begin{equation*}
P_{1}\in \overline{G}_{U_{1},U_{2}}^{+},\text{ }P_{2}\in \overline{G}%
_{U_{2},U_{3}}^{+},\text{ }...,\text{ }P_{n}\in \overline{G}%
_{U_{n},U_{n+1}}^{+}.
\end{equation*}

\smallskip

\textbf{Proof.} (i) ``$\Longrightarrow $'' The matrix $\left(
P_{1}P_{2}\right) _{U_{1}}^{U_{3}}$ is row-allowable because $P_{1}P_{2}\in 
\overline{G}_{U_{1},U_{3}}^{+}.$ It follows that $\forall i\in U_{1},$ $%
\exists j=j\left( i\right) \in U_{3}$ such that 
\begin{equation*}
\left( P_{1}P_{2}\right) _{ij}>0.
\end{equation*}
\noindent Further, we have 
\begin{equation*}
0<\left( P_{1}P_{2}\right) _{ij}=\sum_{k\in \left\langle m_{2}\right\rangle
}\left( P_{1}\right) _{ik}\left( P_{2}\right) _{kj},\text{ }\forall i\in
U_{1}.
\end{equation*}

\noindent It follows that $\forall i\in U_{1},$ $\exists k=k\left( i\right)
\in \left\langle m_{2}\right\rangle ,$ $\exists j=j\left( i\right) \in U_{3}$
such that 
\begin{equation*}
\left( P_{1}\right) _{ik},\text{ }\left( P_{2}\right) _{kj}>0.
\end{equation*}
\noindent Set 
\begin{equation*}
U_{2}=\bigcup_{i\in U_{1}}\left\{ k\left( i\right) \right\} .
\end{equation*}
\noindent Obviously, we have 
\begin{equation*}
P_{1}\in \overline{G}_{U_{1},U_{2}}^{+}\text{ and }P_{2}\in \overline{G}%
_{U_{2},U_{3}}^{+}.
\end{equation*}

``$\Longleftarrow $'' See [28, Theorem 2.6].

\smallskip

(ii) Induction. 
\endproof%

\smallskip

Theorem 2.8 from [28] is a powerful result for the sum-positive matrices ---
two sections of applications in [28]. This theorem/result is a if-then
theorem --- below we give a better theorem, a if-and-only-if theorem.

\smallskip

\textbf{THEOREM 1.7.} \textit{Let }$P_{1}\in N_{m_{1},m_{2}},$\textit{\ }$%
P_{2}\in N_{m_{2},m_{3}},$ $...,$ $P_{n}\in N_{m_{n},m_{n+1}}.$\textit{\ Let 
}$j\in \left\langle m_{n+1}\right\rangle .$\textit{\ Let }$%
U_{1}=\left\langle m_{1}\right\rangle $\textit{\ and }$U_{n+1}=\left\{
j\right\} .$\textit{\ Then} 
\begin{equation*}
\left( P_{1}P_{2}...P_{n}\right) ^{\left\{ j\right\} }>0
\end{equation*}

\noindent (\textit{i.e., the column }$j$\textit{\ of }$P_{1}P_{2}...P_{n}$%
\textit{\ is positive, and, in this case, }$P_{1}P_{2}...P_{n}$\textit{\ is
a Markov matrix}) \textit{if and only if} 
\begin{equation*}
\exists U_{2},U_{3},...,U_{n},
\end{equation*}
\begin{equation*}
\emptyset \neq U_{2}\subseteq \left\langle m_{2}\right\rangle ,\text{ }%
\emptyset \neq U_{3}\subseteq \left\langle m_{3}\right\rangle ,\text{ }...,%
\text{ }\emptyset \neq U_{n}\subseteq \left\langle m_{n}\right\rangle ,
\end{equation*}
\noindent \textit{such that} 
\begin{equation*}
P_{1}\in \overline{G}_{U_{1},U_{2}}^{+},\text{ }P_{2}\in \overline{G}%
_{U_{2},U_{3}}^{+},\text{ }...,\text{ }P_{n}\in \overline{G}%
_{U_{n},U_{n+1}}^{+}.
\end{equation*}

\smallskip

\textbf{Proof.} Theorem 1.6(ii) 
\begin{equation*}
\text{(}\left( P_{1}P_{2}...P_{n}\right) ^{\left\{ j\right\}
}>0\Longleftrightarrow P_{1}P_{2}...P_{n}\in \overline{G}_{\left\langle
m_{1}\right\rangle ,\left\{ j\right\} }^{+}\text{). 
\endproof
}
\end{equation*}

\smallskip

\textbf{Example 1.8.} Let 
\begin{equation*}
P_{1}=\left( 
\begin{array}{ccc}
1 & 0 & 0 \\ 
0 & 1 & 0 \\ 
0 & 1 & 0%
\end{array}
\right) \text{ and }P_{2}=\left( 
\begin{array}{ccc}
1 & 0 & 0 \\ 
1 & 0 & 0 \\ 
0 & 0 & 1%
\end{array}
\right) .
\end{equation*}

\noindent We have $P_{1}\in G_{\left\langle 3\right\rangle ,\left\langle
2\right\rangle }^{+}$ and $P_{2}\in G_{\left\langle 2\right\rangle ,\left\{
1\right\} }^{+}$ ($P_{1}$ and $P_{2}$ are stochastic matrices). By Theorem
1.7, 
\begin{equation*}
\left( P_{1}P_{2}\right) ^{\left\{ 1\right\} }>0
\end{equation*}
\noindent --- the column 1 of $P_{1}P_{2}$ is positive; so, $P_{1}P_{2}$ is
a Markov matrix (this thing can also be obtained by direct computation).

\smallskip

Let $P\in N_{m,n}.$ Let $\emptyset \neq U\subseteq \left\langle
m\right\rangle $ and $\emptyset \neq V\subseteq \left\langle n\right\rangle $%
. Set (see [28]) 
\begin{equation*}
U\rightarrow V\text{ or, if confusion can arise, }U\overset{P}{\rightarrow }V
\end{equation*}

\noindent if $\forall i\in U,$ $\exists j\in V$ such that $P_{ij}>0.$ Set 
\begin{equation*}
V\leftarrow U\text{ or, if confusion can arise, }V\overset{P}{\leftarrow }U
\end{equation*}

\noindent if $U\rightarrow V$. Obviously, if $U\rightarrow V$, then $%
U\rightarrow W,$ $\forall W,$ $V\subseteq W\subseteq \left\langle
n\right\rangle ,$ and $P$ is sum-positive on $U\times V$, more generally, on 
$U\times W,$ $\forall W,$ $V\subseteq W\subseteq \left\langle n\right\rangle
.$

\smallskip

\textbf{Remark 1.9 }([28]). Let $P_{1}\in N_{m_{1},m_{2}},$ $P_{2}\in
N_{m_{2},m_{3}},$ $...,$ $P_{n}\in N_{m_{n},m_{n+1}}.$\ Let $j\in
\left\langle m_{n+1}\right\rangle .$\ If 
\begin{equation*}
U_{1}=\left\langle m_{1}\right\rangle \rightarrow U_{2}\rightarrow
...\rightarrow U_{n}\rightarrow U_{n+1}=\left\{ j\right\}
\end{equation*}
\noindent (equivalently, if 
\begin{equation*}
U_{n+1}=\left\{ j\right\} \leftarrow U_{n}\leftarrow ...\leftarrow
U_{2}\leftarrow U_{1}=\left\langle m_{1}\right\rangle \text{),}
\end{equation*}
\noindent where, obviously, $\emptyset \neq U_{2}\subseteq \left\langle
m_{2}\right\rangle ,$ $...,$ $\emptyset \neq U_{n}\subseteq \left\langle
m_{n}\right\rangle ,$ then, by Theorem 1.7, 
\begin{equation*}
\left( P_{1}P_{2}...P_{n}\right) ^{\left\{ j\right\} }>0.
\end{equation*}

\noindent (We used 
\begin{equation*}
U_{1}=\left\langle m_{1}\right\rangle \rightarrow U_{2}\rightarrow
...\rightarrow U_{n}\rightarrow U_{n+1}=\left\{ j\right\} ,
\end{equation*}

\noindent not 
\begin{equation*}
U_{1}=\left\langle m_{1}\right\rangle \overset{P_{1}}{\rightarrow }U_{2}%
\overset{P_{2}}{\rightarrow }...\overset{P_{n-1}}{\rightarrow }U_{n}\overset{%
P_{n}}{\rightarrow }U_{n+1}=\left\{ j\right\} ,
\end{equation*}

\noindent because no confusion can arise.)

\smallskip

\textbf{Example 1.10.} Consider the matrices from Example 1.8. We have 
\begin{equation*}
\left\langle 3\right\rangle \rightarrow \left\langle 2\right\rangle
\rightarrow \left\{ 1\right\}
\end{equation*}
\noindent ($\left\langle 3\right\rangle \overset{P_{1}}{\rightarrow }%
\left\langle 2\right\rangle \overset{P_{2}}{\rightarrow }\left\{ 1\right\} $%
), 
\begin{equation*}
\left\langle 2\right\rangle \rightarrow \left\langle 2\right\rangle
\rightarrow \left\{ 1\right\}
\end{equation*}

\noindent (equivalently, $\left\{ 1\right\} \leftarrow \left\langle
2\right\rangle \leftarrow \left\langle 2\right\rangle $), 
\begin{equation*}
\left\langle 3\right\rangle \overset{P_{1}}{\rightarrow }\left\langle
2\right\rangle \overset{P_{1}}{\rightarrow }\left\langle 2\right\rangle ,
\end{equation*}
\begin{equation*}
\left\langle 2\right\rangle \overset{P_{1}}{\rightarrow }\left\langle
2\right\rangle ,\text{ }\left\{ 1\right\} \overset{P_{2}}{\rightarrow }%
\left\langle 3\right\rangle ,\text{ etc.}
\end{equation*}

\smallskip

Definition 1.3 can be nuanced --- below we do this thing.

\smallskip

Set 
\begin{equation*}
\mathbb{R}^{+}=\left\{ b\left| \text{ }b\in \mathbb{R}\text{ and }b>0\right.
\right\} .
\end{equation*}

\noindent $\mathbb{R}^{+}$ is the set of positive real numbers.

\smallskip

\textbf{Definition 1.11.} Let $P\in N_{m,n}$ ($m,$ $n\geq 1$). Let $a\in 
\mathbb{R}^{+}$. Let $\emptyset \neq U\subseteq \left\langle m\right\rangle $
and $\emptyset \neq V\subseteq \left\langle n\right\rangle $. We say that $P$
is $a$-\textit{sum-positive on }$U\times V$ if 
\begin{equation*}
\sum\limits_{j\in V}P_{ij}\geq a,\text{ }\forall i\in U.
\end{equation*}

\smallskip

For the $a$-sum-positive matrices on $U\times V,$ $a\in \mathbb{R}^{+},$ $%
\emptyset \neq U\subseteq \left\langle m\right\rangle $, $\emptyset \neq
V\subseteq \left\langle n\right\rangle ,$ we will use the generic name ``$a$%
-sum-positive matrices''.

\smallskip

Let $a\in \mathbb{R}^{+}$. Let $\emptyset \neq U\subseteq \left\langle
m\right\rangle $ and $\emptyset \neq V\subseteq \left\langle n\right\rangle $%
. Set 
\begin{equation*}
G_{a,U,V}^{+}=G_{a,U,V}^{+}\left( m,n\right) =\left\{ P\left| \text{ }P\in
S_{m,n}\text{ and is }a\text{-sum-positive on }U\times V\right. \right\}
\end{equation*}

\noindent and 
\begin{equation*}
\overline{G}_{a,U,V}^{+}=\overline{G}_{a,U,V}^{+}\left( m,n\right) =\left\{
P\left| \text{ }P\in N_{m,n}\text{ and is }a\text{-sum-positive on }U\times
V\right. \right\}
\end{equation*}

\noindent --- $G_{a,U,V}^{+}$ and $\overline{G}_{a,U,V}^{+}$ are two
important sets of $a$-sum-positive matrices; obviously, 
\begin{equation*}
G_{a,U,V}^{+}\subseteq S_{m,n},\text{ }\overline{G}_{a,U,V}^{+}\subseteq
N_{m,n},
\end{equation*}
\begin{equation*}
G_{a,U,V}^{+}\subseteq \overline{G}_{a,U,V}^{+},
\end{equation*}
\begin{equation*}
\bigcup_{a\in \mathbb{R}^{+}}G_{a,U,V}^{+}=G_{U,V}^{+},\text{ }\bigcup_{a\in 
\mathbb{R}^{+}}\overline{G}_{a,U,V}^{+}=\overline{G}_{U,V}^{+}.
\end{equation*}

\smallskip

\textbf{Remark 1.12.} Let $P\in N_{m,n}$. Let $\emptyset \neq U\subseteq
\left\langle m\right\rangle $ and $\emptyset \neq V\subseteq \left\langle
n\right\rangle $.

\smallskip

(a) Let $a\in \mathbb{R}^{+}$. $P$ is $a$-sum-positive on $U\times
V\Longrightarrow P$ is sum-positive on $U\times V;$ $P$ is sum-positive on $%
U\times V\nRightarrow P$ is $a$-sum-positive on $U\times V.$

\smallskip

(b) (An equivalent formulation, using words, for ``$\bigcup\limits_{a\in 
\mathbb{R}^{+}}\overline{G}_{a,U,V}^{+}=\overline{G}_{U,V}^{+}$''.) $P$ is $%
a $-sum-positive on $U\times V$ for some $a\in \mathbb{R}^{+}%
\Longleftrightarrow P$ is sum-positive on $U\times V.$

\smallskip

(c) Which of the two notions, that of a sum-positive matrix on $U\times V$
and that of an $a$-sum-positive matrix on $U\times V$, is more general? ---
a question for the reader.

\smallskip

The next result is a basic one for the $a$-sum-positive matrices --- it is
somehow similar to Theorem 1.6, ``$\Longleftarrow $''.

\smallskip

\textbf{THEOREM 1.13.} (i) \textit{Let} $P_{1}\in \overline{G}%
_{a_{1},U_{1},U_{2}}^{+}\subseteq N_{m_{1},m_{2}}$ \textit{and }$P_{2}\in 
\overline{G}_{a_{2},U_{2},U_{3}}^{+}\subseteq N_{m_{2},m_{3}}.$\textit{\ Then%
} 
\begin{equation*}
P_{1}P_{2}\in \overline{G}_{a_{1}a_{2},U_{1},U_{3}}^{+}.
\end{equation*}

\smallskip

(ii) (\textit{a generalization of }(i)) \textit{Let }$P_{1}\in \overline{G}%
_{a_{1},U_{1},U_{2}}^{+}\subseteq N_{m_{1},m_{2}},$\textit{\ }$P_{2}\in 
\overline{G}_{a_{2},U_{2},U_{3}}^{+}\subseteq N_{m_{2},m_{3}},$ $...,$ $%
P_{n}\in \overline{G}_{a_{n},U_{n},U_{n+1}}^{+}\subseteq N_{m_{n},m_{n+1}}.$ 
\textit{Then} 
\begin{equation*}
P_{1}P_{2}...P_{n}\in \overline{G}_{a_{1}a_{2}...a_{n},U_{1},U_{n+1}}^{+}.
\end{equation*}

\smallskip

\textbf{Proof.} (i) Let $i\in U_{1}.$ We have 
\begin{equation*}
\sum_{j\in U_{3}}\left( P_{1}P_{2}\right) _{ij}=\sum_{j\in U_{3}}\sum_{k\in
\left\langle m_{2}\right\rangle }\left( P_{1}\right) _{ik}\left(
P_{2}\right) _{kj}\geq \sum_{j\in U_{3}}\sum_{k\in U_{2}}\left( P_{1}\right)
_{ik}\left( P_{2}\right) _{kj}=
\end{equation*}
\begin{equation*}
=\sum_{k\in U_{2}}\left( P_{1}\right) _{ik}\sum_{j\in U_{3}}\left(
P_{2}\right) _{kj}\geq a_{2}\sum_{k\in U_{2}}\left( P_{1}\right) _{ik}\geq
a_{1}a_{2}.
\end{equation*}

\smallskip

(ii) Induction. 
\endproof%

\smallskip

Below we give an important result for the $a$-sum-positive matrices --- it
is somehow similar to Theorem 1.7, ``$\Longleftarrow $''.

\smallskip

\textbf{THEOREM 1.14.} \textit{Let }$P_{1}\in \overline{G}%
_{a_{1},U_{1},U_{2}}^{+}\subseteq N_{m_{1},m_{2}},$\textit{\ }$P_{2}\in 
\overline{G}_{a_{2},U_{2},U_{3}}^{+}\subseteq N_{m_{2},m_{3}},$ $...,$ $%
P_{n}\in \overline{G}_{a_{n},U_{n},U_{n+1}}^{+}\subseteq N_{m_{n},m_{n+1}}.$%
\textit{\ Suppose that }$U_{1}=\left\langle m_{1}\right\rangle $\textit{\
and }$U_{n+1}=\left\{ j\right\} $, \textit{where} $j\in \left\langle
m_{n+1}\right\rangle .$\textit{\ Then} 
\begin{equation*}
\left( P_{1}P_{2}...P_{n}\right) ^{\left\{ j\right\} }\geq \left( 
\begin{array}{c}
a_{1}a_{2}...a_{n} \\ 
a_{1}a_{2}...a_{n} \\ 
\vdots \\ 
a_{1}a_{2}...a_{n}%
\end{array}
\right)
\end{equation*}

\noindent (\textit{each entry of the column }$j$\textit{\ of }$%
P_{1}P_{2}...P_{n}$\textit{\ is greater or equal to }$a_{1}a_{2}...a_{n}$, 
\textit{and we have} 
\begin{equation*}
\mu \left( P_{1}P_{2}...P_{n}\right) \geq a_{1}a_{2}...a_{n}\text{).}
\end{equation*}

\smallskip

\textbf{Proof.} Theorem 1.13(ii). 
\endproof%

\smallskip

Below we give two applications of Theorem 1.14.

\smallskip

\textbf{Application 1.15.} (An application to Markov chains, see, \textit{%
e.g.}, [3, p. 257] and [18, p. 125].) Let $P\in S_{n}$ ($n\geq 1$). Let $y$
be a column vector with $n$ components. Then 
\begin{equation*}
\min_{i\in \left\langle n\right\rangle }\left( Py\right) _{i}\geq \min_{i\in
\left\langle n\right\rangle }y_{i},
\end{equation*}
\begin{equation*}
\max_{i\in \left\langle n\right\rangle }\left( Py\right) _{i}\leq \max_{i\in
\left\langle n\right\rangle }y_{i}
\end{equation*}
\noindent (the column vector $Py$ is or can be more contracted than the
column vector $y$).

The first inequality is obvious when $\min\limits_{i\in \left\langle
n\right\rangle }y_{i}=0.$ Now we consider that $\min\limits_{i\in
\left\langle n\right\rangle }y_{i}>0.$ We have $P\in G_{1,\left\langle
n\right\rangle ,\left\langle n\right\rangle }^{+}.$ Considering the column
matrix 
\begin{equation*}
Y\in N_{n,1,}\text{ }Y_{i1}=y_{i},\text{ }\forall i\in \left\langle
n\right\rangle ,
\end{equation*}
\noindent we have $Y\in \overline{G}_{a_{2},\left\langle n\right\rangle
,\left\langle 1\right\rangle }^{+},$ where 
\begin{equation*}
a_{2}=\min\limits_{i\in \left\langle n\right\rangle }y_{i}.
\end{equation*}
\noindent By Theorem 1.14 we have ($a_{1}=1$) 
\begin{equation*}
PY\in \overline{G}_{a_{2},\left\langle n\right\rangle ,\left\langle
1\right\rangle }^{+}
\end{equation*}
\noindent --- therefore, 
\begin{equation*}
\left( Py\right) _{i}\geq a_{2}=\min\limits_{i\in \left\langle
n\right\rangle }y_{i},\text{ }\forall i\in \left\langle n\right\rangle ,
\end{equation*}
\noindent and, further, 
\begin{equation*}
\min_{i\in \left\langle n\right\rangle }\left( Py\right) _{i}\geq \min_{i\in
\left\langle n\right\rangle }y_{i}.
\end{equation*}

To prove the second inequality, setting 
\begin{equation*}
d=\max_{i\in \left\langle n\right\rangle }y_{i},
\end{equation*}
\begin{equation*}
g=\left( 
\begin{array}{c}
d+z \\ 
d+z \\ 
\vdots \\ 
d+z%
\end{array}
\right) \in \mathbb{R}^{n},
\end{equation*}
\noindent where $z>0$ ($z>0\Longrightarrow g\in \left( \mathbb{R}^{+}\right)
^{n}$), and 
\begin{equation*}
h=g-y\text{ (}h\in \left( \mathbb{R}^{+}\right) ^{n}\text{),}
\end{equation*}
\noindent we have 
\begin{equation*}
Pg=g
\end{equation*}
\noindent and 
\begin{equation*}
h_{i}\geq z,\text{ }\forall i\in \left\langle n\right\rangle .
\end{equation*}
\noindent Considering the column matrix 
\begin{equation*}
H\in N_{n,1},\text{ }H_{i1}=h_{i},\text{ }\forall i\in \left\langle
n\right\rangle ,
\end{equation*}
\noindent we have $H\in \overline{G}_{z,\left\langle n\right\rangle
,\left\langle 1\right\rangle }^{+}$. By Theorem 1.14 we have ($P\in
G_{1,\left\langle n\right\rangle ,\left\langle n\right\rangle }^{+}$) 
\begin{equation*}
PH\in \overline{G}_{z,\left\langle n\right\rangle ,\left\langle
1\right\rangle }^{+}
\end{equation*}
\noindent --- therefore, 
\begin{equation*}
\left( Ph\right) _{i}\geq z,\text{ }\forall i\in \left\langle n\right\rangle
.
\end{equation*}
\noindent Further, 
\begin{equation*}
z\leq \left( Ph\right) _{i}=\left( Pg\right) _{i}-\left( Py\right)
_{i}=g_{i}-\left( Py\right) _{i},\text{ }\forall i\in \left\langle
n\right\rangle .
\end{equation*}
\noindent Further, 
\begin{equation*}
\left( Py\right) _{i}\leq g_{i}-z=d,\text{ }\forall i\in \left\langle
n\right\rangle .
\end{equation*}
\noindent So, 
\begin{equation*}
\max_{i\in \left\langle n\right\rangle }\left( Py\right) _{i}\leq \max_{i\in
\left\langle n\right\rangle }y_{i}.
\end{equation*}

Note that it is easy to prove the two inequalities without using Theorem
1.14 --- we only prove the second one. We have 
\begin{equation*}
\left( Py\right) _{i}=\sum_{j\in \left\langle n\right\rangle
}P_{ij}y_{j}\leq \sum_{j\in \left\langle n\right\rangle }P_{ij}\max_{k\in
\left\langle n\right\rangle }y_{k}=\max_{k\in \left\langle n\right\rangle
}y_{k}\sum_{j\in \left\langle n\right\rangle }P_{ij}=\max_{k\in \left\langle
n\right\rangle }y_{k},\text{ }\forall i\in \left\langle n\right\rangle .
\end{equation*}
\noindent So, 
\begin{equation*}
\max_{i\in \left\langle n\right\rangle }\left( Py\right) _{i}\leq \max_{i\in
\left\langle n\right\rangle }y_{i}.
\end{equation*}
\noindent Note, moreover, that if $P\in S_{n}$ and $Q\in N_{n,p}$ ($n,$ $%
p\geq 1;$ in particular, $Q\in S_{n}$), then, applying the above result to 
\begin{equation*}
P\text{ and }Q^{\left\{ j\right\} },\forall j\in \left\langle p\right\rangle
,
\end{equation*}
we obtain that $PQ$ is or can be more contracted than $Q$ (the rows of $PQ$
are or can be closer than those of $Q$).

\smallskip

\textbf{Application 1.16.} (An application in the reliability theory --- to
give an upper bound for the reliability of a system with independent and
linearly connected components which fails if at least $k$ consecutive
components fail.) Consider the consecutive-k-out-of-$\infty $: F Markov
chain (see, \textit{e.g.}, [25, Definition 2.10]) with state space $%
S=\left\langle \left\langle k\right\rangle \right\rangle $ ($\left\langle
\left\langle k\right\rangle \right\rangle =\left\{ 0,1,...,k\right\} $), $%
k\geq 1,$ and transition matrices $P_{n},$ $n\geq 1,$ where 
\begin{equation*}
P_{n}=\left( 
\begin{array}{cccccc}
p_{n} & q_{n} &  &  &  &  \\ 
p_{n} &  & q_{n} &  &  &  \\ 
. & . & . & . & . & . \\ 
&  &  &  &  &  \\ 
p_{n} &  &  &  &  & q_{n} \\ 
&  &  &  &  & 1%
\end{array}
\right) ,\text{ }\forall n\geq 1
\end{equation*}
\noindent ($p_{n},$ $q_{n}\geq 0,$ $p_{n}+q_{n}=1,$ $\forall n\geq 1$; $%
P_{n}\in S_{k+1},$ $\forall n\geq 1$). Taking the matrices 
\begin{equation*}
P_{m+1},\text{ }P_{m+2},\text{ }...,\text{ }P_{m+k}\text{ (}m\geq 0\text{)}
\end{equation*}
\noindent and sets 
\begin{equation*}
U_{1}=\left\langle \left\langle k\right\rangle \right\rangle ,\text{ }%
U_{2}=\left\langle k\right\rangle ,\text{ }U_{3}=\left\{ 2,3,...,k\right\}
,...,
\end{equation*}
\begin{equation*}
U_{k}=\left\{ k-1,k\right\} ,\text{ }U_{k+1}=\left\{ k\right\} ,
\end{equation*}

\noindent we have 
\begin{equation*}
P_{m+1}\in G_{q_{m+1},U_{1},U_{2}}^{+}\subseteq S_{k+1},\text{ }P_{m+2}\in
G_{q_{m+2},U_{2},U_{3}}^{+}\subseteq S_{k+1},\text{ }...,
\end{equation*}
\begin{equation*}
P_{m+k}\in G_{q_{m+k},U_{k},U_{k+1}}^{+}\subseteq S_{k+1}.
\end{equation*}
\noindent Further, by Theorem 1.14, 
\begin{equation*}
\mu \left( P_{m+1}P_{m+2}...P_{m+k}\right) \geq q_{m+1}q_{m+2}...q_{m+k}.
\end{equation*}
\noindent By Theorem 1.13 from [22] we have 
\begin{equation*}
\alpha \left( P_{m+1}P_{m+2}...P_{m+k}\right) =\mu \left(
P_{m+1}P_{m+2}...P_{m+k}\right) \geq q_{m+1}q_{m+2}...q_{m+k},
\end{equation*}
\noindent where $\alpha $ is the Dobrushin ergodicity coefficient --- if $%
Q\in S_{n,p},$%
\begin{equation*}
\alpha \left( Q\right) =\min_{1\leq i,j\leq n}\sum_{b=1}^{p}\min \left(
Q_{ib},Q_{jb}\right)
\end{equation*}

\noindent (see, \textit{e.g.}, [18, pp. 56--57] and [22]). We obtained, here
using the $G^{+}$ method, the same inequality for $\alpha \left(
P_{m+1}P_{m+2}...P_{m+k}\right) $ as that from the proof of Theorem 2.11
from [25], there using the $G$ method, more precisely, Theorem 1.8 from [24]
--- further, see, in [25], Theorem 2.11 and its proof and other things. (For
the $G$ method, see [24], see also [23] and [29].)

\smallskip

To give other applications of Theorem 1.14, \textit{e.g.}, to give
applications in the reliability theory and to waiting time random variables,
see, \textit{e.g.}, [22] and [25]-[27].

\smallskip

Set 
\begin{equation*}
\text{Par}\left( E\right) =\left\{ \Delta \left| \text{ }\Delta \text{ is a
partition of }E\right. \right\} ,
\end{equation*}

\noindent where $E$ is a nonempty set. We will agree that the partitions do
not contain the empty set. $\left( E\right) $ is the improper (degenerate)
partition of $E$, is the coarsest partition of $E$ --- we need this
partition.

\smallskip

\textbf{Definition\ 1.17.} Let $\Delta _{1},\Delta _{2}\in $Par$\left(
E\right) .$ We say that $\Delta _{1}$\textit{\ is finer than }$\Delta _{2}$
if $\forall V\in \Delta _{1},$ $\exists W\in \Delta _{2}$ such that $%
V\subseteq W.$

\smallskip

Write $\Delta _{1}\preceq \Delta _{2}$ when $\Delta _{1}$ is finer than $%
\Delta _{2}.$

\smallskip

Let $W$ be a finite nonempty set. Suppose that $W=\left\{
s_{1},s_{2},...,s_{t}\right\} $. Set 
\begin{equation*}
\left( \left\{ i\right\} \right) _{i\in W}\in \text{Par}\left( W\right) ,%
\text{ }\left( \left\{ i\right\} \right) _{i\in W}=\left( \left\{
s_{1}\right\} ,\left\{ s_{2}\right\} ,...,\left\{ s_{t}\right\} \right) .
\end{equation*}

\noindent $\left( \left\{ i\right\} \right) _{i\in W}$ is the finest
partition of $W.$ \textit{E.g.}, 
\begin{equation*}
\left( \left\{ i\right\} \right) _{i\in \left\langle n\right\rangle }=\left(
\left\{ 1\right\} ,\left\{ 2\right\} ,...,\left\{ n\right\} \right) .
\end{equation*}

\smallskip

\textbf{Definition 1.18} ([28]). Let $P\in N_{m,n}$. Let $\Delta \in $Par$%
\left( \left\langle m\right\rangle \right) $ and $\Sigma \in $Par$\left(
\left\langle n\right\rangle \right) $. We say that $P$ is a $\left[ \Delta %
\right] $\textit{-positive matrix on }$\Sigma $ if $\forall V\in \Sigma ,$ $%
\exists U\in \Delta $ such that $P$ is sum-positive on $U\times V$
(equivalently, if $\forall V\in \Sigma ,$ $\exists U\in \Delta $ such that $%
P_{U}^{V}$ is a row-allowable matrix). A $\left[ \Delta \right] $-positive
matrix on $\left( \left\{ j\right\} \right) _{j\in \left\langle
n\right\rangle }$ is called $\left[ \Delta \right] $\textit{-positive }for
short.

\smallskip

\textbf{Definition 1.19} ([28]). Let $P\in N_{m,n}$. Let $\Delta \in $Par$%
\left( \left\langle m\right\rangle \right) $ and $\Sigma \in $Par$\left(
\left\langle n\right\rangle \right) $. We say that $P$ is a $\Delta $\textit{%
-positive matrix on }$\Sigma $ if $\Delta $ is the least fine partition for
which $P$ is a $\left[ \Delta \right] $-positive matrix on $\Sigma $. A $%
\Delta $-positive matrix on $\left( \left\{ j\right\} \right) _{j\in
\left\langle n\right\rangle }$ is called $\Delta $\textit{-positive} while a 
$\left( \left\langle m\right\rangle \right) $-positive matrix on $\Sigma $
is called \textit{positive\ on }$\Sigma $ for short. A positive matrix on $%
\left( \left\{ j\right\} \right) _{j\in \left\langle n\right\rangle }$ is
called \textit{positive} for short --- in this case, $P\in \overline{G}%
_{\left( \left\langle m\right\rangle \right) ,\left( \left\{ j\right\}
\right) _{j\in \left\langle n\right\rangle }}^{+},$ so $P>0,$ and,
therefore, ``\textit{positive}'' is justified; we have more, namely, $P>0$
if and only if $P\in \overline{G}_{\left( \left\langle m\right\rangle
\right) ,\left( \left\{ j\right\} \right) _{j\in \left\langle n\right\rangle
}}^{+}$.

\smallskip

For the $\left[ \Delta \right] $- and $\Delta $-positive matrices on $\Sigma
,$ $\Delta \in $Par$\left( \left\langle m\right\rangle \right) $, $\Sigma
\in $Par$\left( \left\langle n\right\rangle \right) $, we will use the
generic name ``generalized positive matrices'' --- in [28] these matrices
were differently named...

\smallskip

Let $\Delta \in $Par$\left( \left\langle m\right\rangle \right) $ and $%
\Sigma \in $Par$\left( \left\langle n\right\rangle \right) $. Set (see [28]) 
\begin{equation*}
G_{\Delta ,\Sigma }^{+}=G_{\Delta ,\Sigma }^{+}\left( m,n\right) =\left\{
P\left| \text{ }P\in S_{m,n}\text{ and is a }\left[ \Delta \right] \text{%
-positive matrix on }\Sigma \right. \right\}
\end{equation*}

\noindent and 
\begin{equation*}
\overline{G}_{\Delta ,\Sigma }^{+}=\overline{G}_{\Delta ,\Sigma }^{+}\left(
m,n\right) =\left\{ P\left| \text{ }P\in N_{m,n}\text{ and is a }\left[
\Delta \right] \text{-positive matrix on }\Sigma \right. \right\}
\end{equation*}

\noindent --- $G_{\Delta ,\Sigma }^{+}$ and $\overline{G}_{\Delta ,\Sigma
}^{+}$ are two important sets of generalized positive matrices; obviously, $%
G_{\Delta ,\Sigma }^{+}\subseteq S_{m,n},$ $\overline{G}_{\Delta ,\Sigma
}^{+}\subseteq N_{m,n},$ and $G_{\Delta ,\Sigma }^{+}\subseteq \overline{G}%
_{\Delta ,\Sigma }^{+}.$

\smallskip

When we study or even when we construct products of nonnegative matrices (in
particular, products of stochastic matrices) using $G_{U,V}^{+},$ $\overline{%
G}_{U,V}^{+},$ $G_{a,U,V}^{+},$ $\overline{G}_{a,U,V}^{+},$ $G_{\Delta
,\Sigma }^{+},$ or $\overline{G}_{\Delta ,\Sigma }^{+}$, we will refer this
as the $G^{+}$\textit{\ method}. $G^{+}$ comes from the verb \textit{to group%
} and its derivatives and the adjective \textit{positive}.

\smallskip

The next result is a basic one for the generalized positive matrices, and is
somehow similar to Theorems 1.6, ``$\Longleftarrow $'', and 1.13.

\smallskip

\textbf{THEOREM 1.20} ([28]). (i) \textit{Let} $P_{1}\in \overline{G}%
_{\Delta _{1},\Delta _{2}}^{+}\subseteq N_{m_{1},m_{2}}$ \textit{and }$%
P_{2}\in \overline{G}_{\Delta _{2},\Delta _{3}}^{+}\subseteq
N_{m_{2},m_{3}}. $\textit{\ Then} 
\begin{equation*}
P_{1}P_{2}\in \overline{G}_{\Delta _{1},\Delta _{3}}^{+}\subseteq
N_{m_{1},m_{3}}.
\end{equation*}

\smallskip

(ii) (\textit{a generalization of }(i)) \textit{Let }$P_{1}\in \overline{G}%
_{\Delta _{1},\Delta _{2}}^{+}\subseteq N_{m_{1},m_{2}},$\textit{\ }$%
P_{2}\in \overline{G}_{\Delta _{2},\Delta _{3}}^{+}\subseteq
N_{m_{2},m_{3}}, $ $...,$ $P_{n}\in \overline{G}_{\Delta _{n},\Delta
_{n+1}}^{+}\subseteq N_{m_{n},m_{n+1}}.$\textit{\ Then} 
\begin{equation*}
P_{1}P_{2}...P_{n}\in \overline{G}_{\Delta _{1},\Delta _{n+1}}^{+}\subseteq
N_{m_{1},m_{n+1}}.
\end{equation*}

\smallskip

\textbf{Proof.} See [28, Theorem 2.15]. 
\endproof%

\smallskip

The next result is somehow similar to Theorems 1.7, ``$\Longleftarrow $'',
and 1.14.

\smallskip

\textbf{THEOREM\ 1.21} ([28]). \textit{Let }$P_{1}\in \overline{G}_{\Delta
_{1},\Delta _{2}}^{+}\subseteq N_{m_{1},m_{2}},$\textit{\ }$P_{2}\in 
\overline{G}_{\Delta _{2},\Delta _{3}}^{+}\subseteq N_{m_{2},m_{3}},$ $...,$ 
$P_{n}\in \overline{G}_{\Delta _{n},\Delta _{n+1}}^{+}\subseteq
N_{m_{n},m_{n+1}}.$\textit{\ Suppose that }$\Delta _{1}=\left( \left\langle
m_{1}\right\rangle \right) $\textit{\ and }$\Delta _{n+1}=\left( \left\{
j\right\} \right) _{j\in \left\langle m_{n+1}\right\rangle }.$\textit{\ Then}
\begin{equation*}
P_{1}P_{2}...P_{n}>0.
\end{equation*}

\smallskip

\textbf{Proof.} See [28, Theorem 2.17]). 
\endproof%

\smallskip

\textbf{Remark 1.22.} Theorems 1.20 and 1.21 are if-then theorems/results
--- the conditions from/of these theorems are only sufficient... Indeed, if 
\begin{equation*}
P_{1}=\left( 
\begin{array}{cc}
1 & 0 \\ 
1 & 0%
\end{array}
\right) \text{ and }P_{2}=\left( 
\begin{array}{cc}
1 & 1 \\ 
0 & 0%
\end{array}
\right) ,
\end{equation*}
\noindent then $P_{1}P_{2}$ is positive, 
\begin{equation*}
P_{1}P_{2}=\left( 
\begin{array}{cc}
1 & 1 \\ 
1 & 1%
\end{array}
\right) ,
\end{equation*}
\noindent but $\nexists \Delta _{2}\in $Par$\left( \left\langle
2\right\rangle \right) $ such that $P_{1}\in \overline{G}_{\Delta
_{1},\Delta _{2}}^{+}$ and $P_{2}\in \overline{G}_{\Delta _{2},\Delta
_{3}}^{+},$ where $\Delta _{1}=\left( \left\langle 2\right\rangle \right) $
and $\Delta _{3}=\left( \left\{ 1\right\} ,\left\{ 2\right\} \right) $ ---
for $\Delta _{2}=\left( \left\langle 2\right\rangle \right) ,$ we have $%
P_{1}\in \overline{G}_{\Delta _{1},\Delta _{2}}^{+}$ and $P_{2}\notin 
\overline{G}_{\Delta _{2},\Delta _{3}}^{+}$ while for $\Delta _{2}=\left(
\left\{ 1\right\} ,\left\{ 2\right\} \right) ,$ we have $P_{1}\notin 
\overline{G}_{\Delta _{1},\Delta _{2}}^{+}$ and $P_{2}\in \overline{G}%
_{\Delta _{2},\Delta _{3}}^{+}.$ Note that $P_{1}$ is row-allowable while $%
P_{2}$ is not. When all matrices are row-allowable, it is possible as
conditions from Theorems 1.20 and 1.21 be necessary and sufficient... --- a
problem.

\smallskip

Below we give a few combinatorial problems on the $G^{+}$ method.

\smallskip

\textbf{Problem 1.23.} Fix $m,$ $n\geq 1.$ Fix $P\in N_{m,n}.$ Find the
number of pairs $\left( U,V\right) ,$ $\emptyset \neq U\subseteq
\left\langle m\right\rangle ,$ $\emptyset \neq V\subseteq \left\langle
n\right\rangle ,$ such that $P$ is sum-positive on $U\times V$
(equivalently, find the number of row-allowable submatrices of $P$).

\smallskip

\textbf{Problem 1.24.} Fix $a>0.$ Fix $m,$ $n\geq 1.$ Fix $P\in N_{m,n}.$
Find the number of pairs $\left( U,V\right) ,$ $\emptyset \neq U\subseteq
\left\langle m\right\rangle ,$ $\emptyset \neq V\subseteq \left\langle
n\right\rangle ,$ such that $P$ is $a$-sum-positive on $U\times V$
(equivalently, find the number of row-allowable submatrices of $P$ with row
sums greater or equal to $a$).

\smallskip

\textbf{Problem 1.25.} Fix $m,$ $n\geq 1.$ Fix $P\in N_{m,n}.$ Find the
number of pairs $\left( \Delta ,\Sigma \right) ,$ $\Delta \in $Par$\left(
\left\langle m\right\rangle \right) ,$ $\Sigma \in $Par$\left( \left\langle
n\right\rangle \right) ,$ such that $P$ is $\left[ \Delta \right] $-positive
on $\Sigma .$ Similarly, find the number of pairs $\left( \Delta ,\Sigma
\right) ,$ $\Delta \in $Par$\left( \left\langle m\right\rangle \right) ,$ $%
\Sigma \in $Par$\left( \left\langle n\right\rangle \right) ,$ such that $P$
is $\Delta $-positive on $\Sigma .$

\smallskip

\textbf{Remark 1.26}. The number of partitions of a set with $n$ elements is
given by the Bell number $B_{n},$ see, \textit{e.g.}, [20, pp. 418--419].

\smallskip

The $G^{+}$ method has many applications. Besides those from [28] --- two
sections of applications (to Markov chains...) --- and, in this article,
Applications 1.15 and 1.16 and those from Section 2 to deterministic finite
automata, this method can be applied in other places, such as, to see if a
nonempty set of stochastic square matrices having the same dimension is a
consensus set (see, \textit{e.g.}, [6] for the consensus sets...) and to see
if a nonempty set of nonnegative square matrices having the same dimension
is a primitive set (see, \textit{e.g.}, [7], [16], [31], and [41] for the
primitive sets...) --- the primitive sets can be useful, \textit{e.g.}, for
the study of Lyapunov exponents, see, \textit{e.g.}, [31, pp. 764--765] and
[32, pp. 4457...].

\bigskip

\begin{center}
\textbf{2. \v{C}ERN\'{Y} CONJECTURE}
\end{center}

\bigskip

In this section, we consider three combinatorial conjectures, the \v{C}ern%
\'{y} conjecture --- a (combinatorial) conjecture for the synchronizable
deterministic finite automata ---, the extended \v{C}ern\'{y} conjecture,
and the more extended \v{C}ern\'{y} conjecture. To prove the \v{C}ern\'{y}
conjecture for classes of automata, several methods are used, two of them
being the compression method and extension method --- these two are
submethods of the $G^{+}$ method, and, to illustrate this thing, we give two
examples.

\smallskip

\textbf{Definition 2.1.} (A suitable definition of a deterministic finite
automaton (without an initial state and a set of accepting/final states) for
our interest, for the \v{C}ern\'{y} conjecture, see, \textit{e.g.}, [16] ---
compare with [33, Definition 1].) Let $Q$ and $\Sigma $ be two finite
nonempty sets. Let $\delta :Q\times \Sigma \longrightarrow Q$ be a function.
The triple $\mathcal{A}=\left( Q,\Sigma ,\delta \right) $ is called a 
\textit{deterministic finite automaton} --- further, in general, for
simplification, we will omit ``deterministic finite''. $Q$ is called the 
\textit{state set} (\textit{of automaton} $\mathcal{A}$), $\Sigma $ is
called the \textit{alphabet}, and $\delta $ is called the \textit{transition
function}.

\smallskip

Each automaton, by means of words, evolves in discrete time. A word of an
automaton is a concatenation of letters of the automaton. Considering an
automaton $\mathcal{A}=\left( Q,\Sigma ,\delta \right) $ and, \textit{e.g.},
the (input) word $w=a_{1}a_{2}a_{3}$ ($a_{1},$ $a_{2},$ $a_{3}\in \Sigma $),
and supposing that the automaton is in the state, say, $q$ ($q\in Q$), the
automaton passes from the state $q$ to the state $\delta \left(
q,a_{1}\right) =:q_{1}$ ($q_{1}\in Q$), then from the state $q_{1}$ to the
state $\delta \left( q_{1},a_{2}\right) =:q_{2},$ then from the state $q_{2}$
to the state $\delta \left( q_{2},a_{3}\right) =:q_{3}.$

\smallskip

Let $\mathcal{A}=\left( Q,\Sigma ,\delta \right) $ be an automaton. Let $%
q\in Q.$ Let $w=a_{1}a_{2}...a_{u}$ be a word of the automaton ($u\geq 1$).
Set (see, \textit{e.g.}, [39]-[40]...) 
\begin{equation*}
q.w=\delta \left( ...\delta \left( \delta \left( q,a_{1}\right)
,a_{2}\right) ...,a_{u}\right) .
\end{equation*}
\noindent $q.w$ is (called) the \textit{image of} (\textit{state}) $q$ 
\textit{under the action of }(\textit{word}) $w$ (see, \textit{e.g.}, [16]). 
\textit{E.g.}, considering the example from the previous paragraph, we have 
\begin{equation*}
q.w=q_{3}
\end{equation*}
\noindent because 
\begin{equation*}
q.w=\delta \left( \delta \left( \delta \left( q,a_{1}\right) ,a_{2}\right)
,a_{3}\right) =\delta \left( \delta \left( q_{1},a_{2}\right) ,a_{3}\right)
=\delta \left( q_{2},a_{3}\right) =q_{3}.
\end{equation*}

\smallskip

\textbf{Definition 2.2.} (See, \textit{e.g.}, [37]-[40] and [42].) Let $%
\mathcal{A}=\left( Q,\Sigma ,\delta \right) $ be an automaton. The automaton 
$\mathcal{A}$ is called \textit{synchronizable} or \textit{synchronizing} if
there exist a word $w$ and a state $s\in Q$ such that $q.w=s,$ $\forall q\in
Q$ --- $w$ is called a \textit{reset} or a \textit{synchronizing word} and $%
s $ is called a \textit{sink}.

\smallskip

\textbf{Definition 2.3.} (See, \textit{e.g.}, [28, Definition 3.17].) Let $%
B\in N_{m,n}$. We say that $B$ is a $\left( 0,1\right) $\textit{-matrix} if $%
B_{ij}\in \left\{ 0,1\right\} ,$ $\forall i\in \left\langle m\right\rangle ,$
$\forall j\in \left\langle n\right\rangle .$

\smallskip

Set 
\begin{equation*}
S_{m,n}^{\left( 0,1\right) }=\left\{ P\left| P\in S_{m,n}\text{ and is a }%
\left( 0,1\right) \text{-matrix}\right. \right\}
\end{equation*}
\noindent (for $S_{m,n}$, see Section 1; $S_{m,n}^{\left( 0,1\right) }$ is
the set of stochastic $m\times n$ $\left( 0,1\right) $-matrices) and 
\begin{equation*}
S_{n}^{\left( 0,1\right) }=S_{n,n}^{\left( 0,1\right) }.
\end{equation*}
\noindent We have 
\begin{equation*}
\left| S_{m,n}^{\left( 0,1\right) }\right| =n^{m}\text{ and }\left|
S_{n}^{\left( 0,1\right) }\right| =n^{n}
\end{equation*}
\noindent ($\left| \cdot \right| $ is the cardinality). \textit{E.g.}, 
\begin{equation*}
\left| S_{100}^{\left( 0,1\right) }\right| =100^{100}=10^{200}
\end{equation*}
\noindent --- a big number!

\smallskip

Further, we consider automata with state set $\left\langle n\right\rangle $
and alphabet $\left\langle t\right\rangle .$

\smallskip

Let $\mathcal{A}=\left( \left\langle n\right\rangle ,\left\langle
t\right\rangle ,\delta \right) $ be an automaton. To each letter $a\in
\left\langle t\right\rangle $ we associate (see, \textit{e.g.}, [5]) a
stochastic $n\times n$ $\left( 0,1\right) $-matrix, $P_{a},$%
\begin{equation*}
\left( P_{a}\right) _{ij}=\left\{ 
\begin{array}{ll}
0\smallskip & \text{if }i.a\neq j, \\ 
1 & \text{if }i.a=j,%
\end{array}
\right.
\end{equation*}
\noindent $\forall i,$ $j\in \left\langle n\right\rangle .$ Obviously, if $%
w=a_{1}a_{2}...a_{u}$ is a word of this automaton, then 
\begin{equation*}
\left( P_{a_{1}}P_{a_{2}}...P_{a_{u}}\right) _{ij}=\left\{ 
\begin{array}{ll}
0\smallskip & \text{if }i.w\neq j, \\ 
1 & \text{if }i.w=j,%
\end{array}
\right.
\end{equation*}
\noindent $\forall i,$ $j\in \left\langle n\right\rangle .$ Due this
association, we can work with stochastic $n\times n$ $\left( 0,1\right) $%
-matrices instead of letters and with products of stochastic $n\times n$ $%
\left( 0,1\right) $-matrices instead of words --- we can work with the
matrices $P_{a},$ $a\in \left\langle t\right\rangle $ ($P_{a}\in
S_{n}^{\left( 0,1\right) },$ $\forall a\in \left\langle t\right\rangle $),
instead of automaton $\mathcal{A}$.

\smallskip

\textbf{Definition 2.4.} Let $\mathcal{C}$ be a (finite or not) nonempty set
of nonnegative matrices. A \textit{product} of matrices from $\mathcal{C}$
is either a matrix from $\mathcal{C}$, and, in this case, the product is
called \textit{degenerate} or \textit{improper}, or a product of two or more
matrices from $\mathcal{C}$, and, in this case, the product is called 
\textit{nondegenerate} or \textit{proper}. A product of matrices from $%
\mathcal{C}$ is called \textit{Markov} or \textit{synchronizing} if it is a
Markov matrix (for the Markov matrices, see Definition 1.4).

\smallskip

\textbf{Definition 2.5.} Let $\emptyset \neq \mathcal{C}\subseteq N_{n}.$
Suppose that $\mathcal{C}$ is a finite set (this condition can be
removed...). We say that $\mathcal{C}$ is a \textit{Markov,} or a \textit{%
synchronizable,} or a \textit{synchronizing set} if there exists a
Markov/synchronizing product of matrices from $\mathcal{C}$.

\smallskip

For the above automaton $\mathcal{A}=\left( \left\langle n\right\rangle
,\left\langle t\right\rangle ,\delta \right) $, setting 
\begin{equation*}
\mathcal{C}=\left\{ P_{a}\left| \text{ }a\in \left\langle t\right\rangle
\right. \right\} \text{ (}\emptyset \neq \mathcal{C}\subseteq S_{n}^{\left(
0,1\right) }\text{),}
\end{equation*}
\noindent we have the following correspondences: 
\begin{equation*}
\text{(automaton) }\mathcal{A}\text{ }\longleftrightarrow \text{ (set) }%
\mathcal{C}\text{;}
\end{equation*}
\begin{equation*}
\text{(state set) }\left\langle n\right\rangle \text{ }\longleftrightarrow 
\text{ index set of matrices from }\mathcal{C}
\end{equation*}
\noindent (for the index sets, see Definition 1.1); 
\begin{equation*}
\text{(alphabet) }\left\langle t\right\rangle \text{ }\longleftrightarrow 
\text{ set of the subscripts of matrices from }\mathcal{C}\text{;}
\end{equation*}
\begin{equation*}
\text{word of }\mathcal{A}\text{ }\longleftrightarrow \text{ product of
matrices from }\mathcal{C}\text{;}
\end{equation*}
\begin{equation*}
\text{reset word (of letters from }\left\langle t\right\rangle \text{) }%
\longleftrightarrow \text{ Markov product (of matrices from }\mathcal{C}%
\text{);}
\end{equation*}
\begin{equation*}
\text{synchronizable automaton }\mathcal{A}\text{, when }\mathcal{A}\text{
is synchronizable }\longleftrightarrow
\end{equation*}
\begin{equation*}
\text{Markov set }\mathcal{C}\text{, when }\mathcal{C}\text{ is a Markov set.%
}
\end{equation*}

\smallskip

Below we give the \v{C}ern\'{y} conjecture, a combinatorial conjecture for
the synchronizable deterministic finite automata.

\smallskip

\textbf{Conjecture 2.6.} (\v{C}ern\'{y} Conjecture [9].) Let $\mathcal{A}%
=\left( \left\langle n\right\rangle ,\left\langle t\right\rangle ,\delta
\right) $ be an automaton, $n\geq 2$ ($t\geq 1$). If $\mathcal{A}$ is
synchronizable, then there exists a reset/synchronizing word of length at
most $\left( n-1\right) ^{2}$ --- equivalently, if $\mathcal{C}$, $\mathcal{C%
}$ corresponds to $\mathcal{A}$, is a Markov set, then there exists a Markov
product of matrices from $\mathcal{C}$ having the length at most $\left(
n-1\right) ^{2}.$

\smallskip

We will extend the \v{C}ern\'{y} conjecture for the Markov sets of
nonnegative matrices --- this idea is not new, see [10].

\smallskip

\textbf{Definition 2.7.} (See, \textit{e.g.}, [18, p. 222] and [35, p. 3].)
Let $P\in N_{m,n}$. Set 
\begin{equation*}
\overline{P}=\left( \overline{P}_{ij}\right) \in N_{m,n},\text{ }\overline{P}%
_{ij}=\left\{ 
\begin{array}{l}
1\text{ if }P_{ij}>0, \\ 
0\text{ if }P_{ij}=0,%
\end{array}
\right.
\end{equation*}

\noindent $\forall i\in \left\langle m\right\rangle ,$ $\forall j\in
\left\langle n\right\rangle .$ We call $\overline{P}$ the \textit{incidence\
matrix of }$P$.

\smallskip

Set 
\begin{equation*}
N_{m,n}^{\left( 0,1\right) }=\left\{ P\left| P\in N_{m,n}\text{ and is a }%
\left( 0,1\right) \text{-matrix}\right. \right\}
\end{equation*}
\noindent (for $N_{m,n},$ see Section 1; $N_{m,n}^{\left( 0,1\right) }$ is
the set of nonnegative $m\times n$ $\left( 0,1\right) $-matrices) and 
\begin{equation*}
N_{n}^{\left( 0,1\right) }=N_{n,n}^{\left( 0,1\right) }.
\end{equation*}

\smallskip

\textbf{Definition 2.8.} (See, \textit{e.g.}, [17, p. 61] and [28,
Definition 3.17].) Let $P\in N_{m,n}$. Let $B\in N_{m,n}^{\left( 0,1\right)
}.$ We say that $P$ \textit{has the pattern} $B$ if $\overline{P}\geq B.$

\smallskip

Set ($m,$ $n\geq 1$) 
\begin{equation*}
\mathcal{R}_{m,n}=\left\{ P\left| P\in N_{m,n}\text{ and is a row-allowable
matrix}\right. \right\}
\end{equation*}
\noindent --- obviously, $S_{m,n}\subseteq \mathcal{R}_{m,n}\subseteq
N_{m,n} $ ---, 
\begin{equation*}
\mathcal{R}_{n}=\mathcal{R}_{n,n},
\end{equation*}
\begin{equation*}
\mathcal{P}_{m,n}\left( B\right) =\left\{ P\left| P\in N_{m,n}\text{ and has
the pattern }B\right. \right\} ,\text{ }\forall B\in S_{m,n}^{\left(
0,1\right) },
\end{equation*}
\begin{equation*}
\mathcal{P}_{n}\left( B\right) =\mathcal{P}_{n,n}\left( B\right) ,\text{ }%
\forall B\in S_{n}^{\left( 0,1\right) },
\end{equation*}
\begin{equation*}
\mathcal{P}_{m,n}=\bigcup_{B\in S_{m,n}^{\left( 0,1\right) }}\mathcal{P}%
_{m,n}\left( B\right) ,
\end{equation*}
\begin{equation*}
\mathcal{P}_{n}=\mathcal{P}_{n,n}.
\end{equation*}

\smallskip

\textbf{THEOREM 2.9.} 
\begin{equation*}
\mathcal{P}_{m,n}=\mathcal{R}_{m,n},\text{ }\forall m,\text{ }n\geq 1
\end{equation*}
\noindent --- \textit{therefore, }$\nexists P\in N_{m,n},$ $P\in \mathcal{R}%
_{m,n}$\textit{\ and }$P\notin \mathcal{P}_{m,n}.$

\smallskip

\textbf{Proof.} ``$\subseteq $'' Obvious (because the stochastic $\left(
0,1\right) $-matrices are row-allowable).

\smallskip

``$\supseteq $'' Let $P\in \mathcal{R}_{m,n}$ ($m,$ $n\geq 1$). $P$ is
row-allowable $\Longrightarrow $ $\overline{P}$ is row-allowable $%
\Longrightarrow $ $\exists D\in S_{m,n}^{\left( 0,1\right) }$ such that $%
\overline{P}\geq D$ $\Longrightarrow $ $P$ has the pattern $D$ $%
\Longrightarrow $ $P\in \mathcal{P}_{m,n}\left( D\right) $ $\Longrightarrow $
$P\in \mathcal{P}_{m,n}.$ 
\endproof%

\smallskip

\textbf{Conjecture 2.10.} (Extended \v{C}ern\'{y} Conjecture.) Let $%
\emptyset \neq \mathcal{C}\subseteq \mathcal{R}_{n}$ be a finite set, $n\geq
2.$ If $\mathcal{C}$ is a Markov set, then there exists a Markov product of
matrices from $\mathcal{C}$ having the length at most $\left( n-1\right)
^{2}.$

\smallskip

The deterministic finite automata are special cases of the probabilistic
finite automata. So, the extended \v{C}ern\'{y} conjecture could be useful
in the combinatorial theory of probabilistic finite automata (for the
probabilistic finite automata, see, \textit{e.g.}, [33]).

\smallskip

\textbf{THEOREM\ 2.11.} \textit{The extended \v{C}ern\'{y} conjecture holds
if and only if the \v{C}ern\'{y} conjecture holds}.

\smallskip

\textbf{Proof.} ``$\Longrightarrow $'' Obvious.

\smallskip

``$\Longleftarrow $'' Theorem 2.9 and [10, Lemma 1]. 
\endproof%

\smallskip

\textbf{Remark 2.12.} In Definition 2.5 we have ``$\emptyset \neq \mathcal{C}%
\subseteq N_{n}$'' while in Conjecture 2.10 we have ``$\emptyset \neq 
\mathcal{C}\subseteq \mathcal{R}_{n}$''. An extension for the finite sets
included in $N_{n}$ of the extended \v{C}ern\'{y} conjecture can be taken
into account. Note that if $\mathcal{C}$ is a Markov set, $\emptyset \neq 
\mathcal{C}\subseteq N_{n}$, then $\exists P\in \mathcal{C}$ such that $P$
is row-allowable. Note that if $A\in N_{m,n}$ and $B\in N_{n,p},$ then

\smallskip

1) if $A$ and $B$ are row-allowable, then $AB$ is row-allowable;

\smallskip

2) if $A$ is row-allowable and $B$ is not row-allowable, then $AB$ is
row-allowable or not --- if 
\begin{equation*}
A=\left( 
\begin{array}{cc}
1 & 0 \\ 
2 & 0%
\end{array}
\right) \text{ (row-allowable) and }B=\left( 
\begin{array}{cc}
1 & 0 \\ 
0 & 0%
\end{array}
\right) \text{ (not row-allowable),}
\end{equation*}
\noindent then $AB=A$ (row-allowable) while if 
\begin{equation*}
A=\left( 
\begin{array}{cc}
0 & 1 \\ 
0 & 2%
\end{array}
\right) \text{ (row-allowable) and }B=\left( 
\begin{array}{cc}
1 & 0 \\ 
0 & 0%
\end{array}
\right) \text{ (not row-allowable),}
\end{equation*}
\noindent then $AB=0$ (not row-allowable);

\smallskip

3) if $A$ is not row-allowable, then $AB$ is not row-allowable.

\smallskip

\textbf{Conjecture 2.13.} (More Extended \v{C}ern\'{y} Conjecture.) Let $%
\emptyset \neq \mathcal{C}\subseteq N_{n}$ be a finite set, $n\geq 2.$ If $%
\mathcal{C}$ is a Markov set, then there exists a Markov product of matrices
from $\mathcal{C}$ having the length at most $\left( n-1\right) ^{2}.$

\smallskip

\textbf{Definition 2.14.} (See, \textit{e.g.}, [21, p. 29].) Let $A$ and $B$
be two (nonnegative, or real, or complex, or...) $m\times n$ matrices. We
say that $A$ and $B$ \textit{have the same zero pattern} when $\forall i\in
\left\langle m\right\rangle ,$ $\forall j\in \left\langle n\right\rangle ,$%
\begin{equation*}
A_{ij}=0\Longleftrightarrow B_{ij}=0.
\end{equation*}

\smallskip

Set $A\overset{zp}{\thicksim }B$ when $A$ and $B$ have the same zero pattern.

\smallskip

\textbf{Remark 2.15.} (a) $\overset{zp}{\thicksim }$ is an equivalence
relation.

\smallskip

(b) Let $A_{1},$ $B_{1}\in N_{m_{1},m_{2}},$ $A_{2},$ $B_{2}\in
N_{m_{2},m_{3}},$ $...,$ $A_{t},$ $B_{t}\in N_{m_{t},m_{t+1}},$ $t\geq 2.$
If 
\begin{equation*}
A_{1}\overset{zp}{\thicksim }B_{1},\text{ }A_{2}\overset{zp}{\thicksim }%
B_{2},\text{ }...,\text{ }A_{t}\overset{zp}{\thicksim }B_{t},
\end{equation*}
\noindent then, by induction, 
\begin{equation*}
A_{1}A_{2}...A_{t}\overset{zp}{\thicksim }B_{1}B_{2}...B_{t}.
\end{equation*}

\smallskip

\textbf{Definition 2.16.} Let $A\in N_{m,n}.$ We say that $A$ is a $\left(
0,1\right) $-\textit{row-sum matrix} if 
\begin{equation*}
\sum\limits_{j=1}^{n}A_{ij}\in \left\{ 0,1\right\} ,\text{ }\forall i\in
\left\langle m\right\rangle .
\end{equation*}

\smallskip

Obviously, the stochastic matrices are $\left( 0,1\right) $-row-sum matrices.

\smallskip

\textbf{Example 2.17.} Let $A\in N_{m,n}\setminus S_{m,n}.$ Let $P$ be a
positive $m\times n$ matrix --- a special case: $P_{ij}=1,$ $\forall i\in
\left\langle m\right\rangle ,$ $\forall j\in \left\langle n\right\rangle .$
Set 
\begin{equation*}
T=P\odot A.
\end{equation*}
\noindent $T$ is the Hadamard product of $P$ and $A,$ \textit{i.e.}, $%
T=\left( T_{ij}\right) =\left( P_{ij}A_{ij}\right) \in N_{m,n}$, $%
T_{ij}=P_{ij}A_{ij},$ $\forall i\in \left\langle m\right\rangle ,$ $\forall
j\in \left\langle n\right\rangle ;$ an interpretation: $P_{ij}$ is the
weight of $A_{ij},$ $\forall i\in \left\langle m\right\rangle ,$ $\forall
j\in \left\langle n\right\rangle .$ Set 
\begin{equation*}
Z_{i}=\sum\limits_{j=1}^{n}T_{ij},\text{ }\forall i\in \left\langle
m\right\rangle .
\end{equation*}
\noindent Set 
\begin{equation*}
B=\left( B_{ij}\right) ,\text{ }B_{ij}=\left\{ 
\begin{array}{ll}
0\smallskip & \text{if }Z_{i}=0, \\ 
\frac{T_{ij}}{Z_{i}} & \text{if }Z_{i}\neq 0,%
\end{array}
\right.
\end{equation*}
\noindent $\forall i\in \left\langle m\right\rangle ,$ $\forall j\in
\left\langle n\right\rangle .$ ($B$ was obtained from $T$ by row
normalization). Obviously, $B$ is a $\left( 0,1\right) $-row-sum matrix, $B$
is a stochastic matrix if $A\in \mathcal{R}_{m,n},$ and 
\begin{equation*}
B\overset{zp}{\thicksim }A
\end{equation*}
\noindent (moreover, $B\overset{zp}{\thicksim }T\overset{zp}{\thicksim }A$).

\smallskip

\textbf{Remark 2.18.} Due to Definitions 2.14 and 2.16, Remark 2.15, and
Example 2.17, the extended \v{C}ern\'{y} conjecture can be reduced to the
case when the matrices are stochastic (we can consider (we can work with)
stochastic matrices instead of row-allowable matrices --- this idea is not
new, see, \textit{e.g.}, [11, p. 2]) while the more extended \v{C}ern\'{y}
conjecture can be reduced to the case when the matrices are $\left(
0,1\right) $-row-sum.

\smallskip

The \v{C}ern\'{y} conjecture holds for several classes of automata, see, 
\textit{e.g.}, [39]-[40]. The extended \v{C}ern\'{y} conjecture holds in all
cases in which the \v{C}ern\'{y} conjecture holds. Besides these cases, in
[28, Sections 4 and 5], there exist several cases in which the extended \v{C}%
ern\'{y} conjecture holds --- see also [17].

\smallskip

Let $\mathcal{A}=\left( \left\langle n\right\rangle ,\left\langle
t\right\rangle ,\delta \right) $ be an automaton. Let $\emptyset \neq
B\subseteq \left\langle n\right\rangle .$ Let $a\in \left\langle
t\right\rangle .$ Set (see, \textit{e.g.}, [39]-[40]) 
\begin{equation*}
B.a=\left\{ i.a\left\vert \text{ }i\in B\right. \right\} \text{ and }%
Ba^{-1}=\left\{ i\left\vert \text{ }i\in \left\langle n\right\rangle \text{, 
}i.a\in B\right. \right\} .
\end{equation*}%
\noindent Set 
\begin{equation*}
\mathcal{M\!}=\mathcal{\!}\left\{ B.a\left\vert \text{ }\emptyset \neq
B\subseteq \left\langle n\right\rangle \mathcal{\!}\text{, }a\in
\left\langle t\right\rangle \right. \right\} \text{ and }\mathcal{M}^{-1}%
\mathcal{\!}=\mathcal{\!}\left\{ Ba^{-1}\left\vert \text{ }\emptyset \neq
B\subseteq \left\langle n\right\rangle \mathcal{\!}\text{, }a\in
\left\langle t\right\rangle \right. \right\} .
\end{equation*}

\smallskip

To prove the \v{C}ern\'{y} conjecture for classes of automata, there exist a
few methods, two of them being, see, \textit{e.g.}, [39]-[40], the
compression method and extension method. These two are submethods of the $%
G^{+}$ method: the compression method corresponds to the case when Theorem
1.7, \textquotedblleft $\Longleftarrow $\textquotedblright , is used doing a
head-to-tail construction of sets, the sets being from $\left\{ \left\langle
n\right\rangle \right\} \cup \mathcal{M},$ while the extension method
corresponds to the case when Theorem 1.7, \textquotedblleft $\Longleftarrow $%
\textquotedblright , is used doing a tail-to-head construction of sets, the
sets being from $\mathcal{M}^{-1}\cup \left\{ \left\{ 1\right\} ,\left\{
2\right\} ,...,\left\{ n\right\} \right\} $ (in this latter case, a starting
set having cardinality 1 is necessary --- a problem). Removing the too
restrictive conditions on the sets on the two methods, we obtain a
generalization of the compression method, see Example 2.21 (see also Theorem
5.12 and its proof in [28]) and a generalization of the extension method,
see Example 2.19 (see also Theorem 4.2 and its proof in [28]), which we also
call, for simplification, the \textit{compression method} and the \textit{%
extension method}, respectively --- these generalized methods are also
submethods of the $G^{+}$ method.

\smallskip

To illustrate the $G^{+}$ method, more precisely, Theorem 1.7, ``$%
\Longleftarrow $'', concerning the \v{C}ern\'{y} conjecture, we give two
examples.

\smallskip

\textbf{Example 2.19.} Consider the automaton $\mathcal{A}=\left(
\left\langle n\right\rangle ,\left\langle t\right\rangle ,\delta \right) ,$ $%
n\geq 2,$ $t\geq 1,$ with 
\begin{equation*}
q.a\geq q,\text{ }\forall q\in \left\langle n\right\rangle ,\text{ }\forall
a\in \left\langle t\right\rangle
\end{equation*}
\noindent (similarly, with 
\begin{equation*}
q.a\leq q,\text{ }\forall q\in \left\langle n\right\rangle ,\text{ }\forall
a\in \left\langle t\right\rangle \text{).}
\end{equation*}
\noindent In this case, 
\begin{equation*}
\mathcal{C}=\left\{ P_{a}\left| \text{ }a\in \left\langle t\right\rangle
\right. \right\}
\end{equation*}
\noindent is a set of upper triangular matrices. Obviously, 
\begin{equation*}
\left( P_{a}\right) _{nn}=1,\text{ }\forall a\in \left\langle t\right\rangle
.
\end{equation*}
\noindent Suppose that $\mathcal{C}$ is a Markov set. We show that there
exists a Markov/synchro-

\noindent nizing product of matrices from $\mathcal{C}$ which has the length
at most $n-1.$ If such a product exists, its column $n$ is positive, the
others being zero columns. We use Theorem 1.7, ``$\Longleftarrow $'', doing
a tail-to-head construction of sets (the extension method). $\exists
P_{a_{1}}\in \mathcal{C}$ ($a_{1}\in \left\langle t\right\rangle $) such
that 
\begin{equation*}
\left\{ n-1,n\right\} \overset{P_{a_{1}}}{\longrightarrow }\left\{ n\right\}
\end{equation*}
\noindent (for notation, see the first paragraph after Example 1.8)--- $%
\left\{ n-1,n\right\} $ is or is not equal to the largest set, $U$, with the
property that $U\overset{P_{a_{1}}}{\longrightarrow }\left\{ n\right\} $
because we can have more, namely,

\begin{equation*}
\left\{ n-k,n-k+1,...,n\right\} \overset{P_{a_{1}}}{\longrightarrow }\left\{
n\right\}
\end{equation*}
\noindent for some $k\in \left\langle n-1\right\rangle $; $P_{a_{1}}$ is a
Markov matrix if $k=n-1$. Indeed, 
\begin{equation*}
\left\{ n\right\} \overset{P_{a}}{\longrightarrow }\left\{ n\right\} ,\text{ 
}\forall a\in \left\langle t\right\rangle ,
\end{equation*}
\noindent and if 
\begin{equation*}
\left\{ n-1\right\} \overset{P_{a}}{\nrightarrow }\left\{ n\right\} ,\text{ }%
\forall a\in \left\langle t\right\rangle ,
\end{equation*}
\noindent then 
\begin{equation*}
\left( P_{a}\right) _{n-1\rightarrow n-1}=1,\text{ }\forall a\in
\left\langle t\right\rangle
\end{equation*}
\noindent (for notation, see the first paragraph before Definition 1.1),
and, further, any product of matrices from $\mathcal{C}$ has the entry $%
\left( n-1,n-1\right) $ equal to $1.$ Contradiction (because $\mathcal{C}$
is a set of upper triangular stochastic $\left( 0,1\right) $-matrices and a
Markov set). Further, $\exists P_{a_{2}}\in \mathcal{C}$ such that 
\begin{equation*}
\left\{ n-2,n-1,n\right\} \overset{P_{a_{2}}}{\longrightarrow }\left\{
n-1,n\right\}
\end{equation*}
\noindent --- we can have more, namely, 
\begin{equation*}
\left\{ n-h,n-h+1,...,n\right\} \overset{P_{a_{2}}}{\longrightarrow }\left\{
n-1,n\right\}
\end{equation*}
\noindent for some $h\in \left\langle n-1\right\rangle -\left\{ 1\right\} $;
by Theorem 1.7, ``$\Longleftarrow $'', $P_{a_{2}}P_{a_{1}}$ is a Markov
matrix if $h=n-1$. Indeed, 
\begin{equation*}
\left\{ n-1,n\right\} \overset{P_{a}}{\longrightarrow }\left\{ n-1,n\right\}
,\text{ }\forall a\in \left\langle t\right\rangle ,
\end{equation*}
\noindent and if 
\begin{equation*}
\left\{ n-2\right\} \overset{P_{a}}{\nrightarrow }\left\{ n-1,n\right\} ,%
\text{ }\forall a\in \left\langle t\right\rangle ,
\end{equation*}
\noindent then 
\begin{equation*}
\left( P_{a}\right) _{n-2\rightarrow n-2}=1,\text{ }\forall a\in
\left\langle t\right\rangle ,
\end{equation*}
\noindent and, further, any product of matrices from $\mathcal{C}$ has the
entry $\left( n-2,n-2\right) $ equal to $1.$ Contradiction. Proceeding in
this way, finally, we obtain 
\begin{equation*}
\left\langle n\right\rangle \overset{P_{a_{n-1}}}{\longrightarrow }%
\left\langle n\right\rangle -\left\{ 1\right\} \overset{P_{a_{n-2}}}{%
\longrightarrow }\left\langle n\right\rangle -\left\langle 2\right\rangle 
\overset{P_{a_{n-3}}}{\longrightarrow }...\overset{P_{a_{2}}}{%
\longrightarrow }\left\{ n-1,n\right\} \overset{P_{a_{1}}}{\longrightarrow }%
\left\{ n\right\}
\end{equation*}
\noindent for some $P_{a_{1}},$ $P_{a_{2}},$ $...,$ $P_{a_{n-1}}\in \mathcal{%
C}.$ By Theorem 1.7, ``$\Longleftarrow $'', we have 
\begin{equation*}
\left( P_{a_{n-1}}P_{a_{n-2}}...P_{a_{1}}\right) ^{\left\{ n\right\} }>0
\end{equation*}
\noindent ($\left( P_{a_{n-1}}P_{a_{n-2}}...P_{a_{1}}\right) ^{\left\{
n\right\} }$ is the column $n$ of $P_{a_{n-1}}P_{a_{n-2}}...P_{a_{1}}$).
Therefore, 
\begin{equation*}
P_{a_{n-1}}P_{a_{n-2}}...P_{a_{1}}
\end{equation*}
\noindent is a Markov product of length $n-1.$

\smallskip

\textbf{Remark 2.20.} Let $P\in S_{n}^{\left( 0,1\right) }$ with $P_{ii}=1$
for some $i\in \left\langle n\right\rangle .$ Suppose that for any $j\in
\left\langle n\right\rangle ,$ $j\neq i,$ there exists a (directed) path
from $j$ to $i.$ (See also Theorem 5.10 in [28] --- for the graph theory,
see, \textit{e.g.}, [2], [4], and [8].) Suppose, for simplification, that $%
P_{nn}=1.$ If $P$ is not upper triangular, by permutation of rows and
columns it can be transformed into an upper triangular matrix --- an
exercise for the reader. So, further, we can suppose that $P$ is upper
triangular --- we suppose this. Let $\mathcal{C}=\left\{ P\right\} .$ This
is a special case of that from Example 2.19. So, 
\begin{equation*}
P^{n-1}
\end{equation*}
\noindent is a Markov matrix (with 
\begin{equation*}
\left( P^{n-1}\right) ^{\left\{ n\right\} }=\left( 
\begin{array}{c}
1 \\ 
1 \\ 
\vdots \\ 
1%
\end{array}
\right) \text{ and }\left( P^{n-1}\right) ^{\left\langle n-1\right\rangle
}=\left( 
\begin{array}{cccc}
0 & 0 & \cdots & 0 \\ 
0 & 0 & \cdots & 0 \\ 
\vdots & \vdots &  & \vdots \\ 
0 & 0 & \cdots & 0%
\end{array}
\right) \text{).}
\end{equation*}
\noindent Note that we can consider any finite set $\mathcal{C}$ with $P\in 
\mathcal{C}$ instead of $\mathcal{C}=\left\{ P\right\} $ --- the \v{C}ern%
\'{y} conjecture holds.

\smallskip

Set 
\begin{equation*}
\left\lfloor u,v\right\rfloor =\left\{ u,u+1,...,v\right\} ,\text{ }\forall
u,v\in \mathbb{Z},\text{ }u\leq v
\end{equation*}
\noindent ($\left\lfloor u,u\right\rfloor =\left\{ u\right\} $). We call $u$
and $v$ the \textit{ends of }(\textit{set}) $\left\lfloor u,v\right\rfloor $%
, $u$ the \textit{left end of }$\left\lfloor u,v\right\rfloor $, and $v$ the 
\textit{right end of }$\left\lfloor u,v\right\rfloor $.

\smallskip

\textbf{Example 2.21.} Let $\mathcal{A}=\left( \left\langle n\right\rangle
,\left\langle t\right\rangle ,\delta \right) ,$ $n\geq 2,$ $t\geq 1,$ be a
monotonic automaton, \textit{i.e.}, an automaton for which the function $%
q\longmapsto q.a$ is nondecreasing, $\forall a\in \left\langle
t\right\rangle $ (equivalently, $i.a\leq j.a,$ $\forall i,$ $j\in
\left\langle n\right\rangle $ with $i\leq j,$ $\forall a\in \left\langle
t\right\rangle $) --- for the monotonic automata, see [1], [39, pp.
875--876], and [40]. On the \v{C}ern\'{y} conjecture, the upper bound for
this automaton is known, is $n-1,$ and the method used to obtain this bound
is neither the compression method nor the extension one, see [1]. Note that
the $G^{+}$ method can be used in some special cases, for some special
monotonic automata, on the \v{C}ern\'{y} conjecture --- we will illustrate
this thing.

\smallskip

Suppose that $\mathcal{C}$ is a Markov set, $\mathcal{C}$ corresponds to $%
\mathcal{A}$, and consider two cases for $\mathcal{C}$.

\smallskip

\textit{Case} 1. $\exists P\in \mathcal{C},$ $\exists !i\in \left\langle
n\right\rangle $ ($\exists !=$ there exists a unique) such that $P_{ii}=1.$
In this case, we will use Theorem 1.7, ``$\Longleftarrow $'', doing a
head-to-tail construction of sets (the compression method) --- we will show
that $P^{b}$ is a Markov matrix for some $b\in \left\langle n-1\right\rangle
.$ Since $P\in \mathcal{C}$, we have $P=P_{a}$ for some $a\in \left\langle
t\right\rangle .$

\smallskip

\textit{Subcase} 1.1. $i-1<n-i$ ---equivalently, $i<\frac{n+1}{2}$; $i-1$ is
the (usual) distance from $i$ to $1$ and $n-i$ is the distance from $i$ to $%
n $.

\smallskip

We show that 
\begin{equation*}
\left\langle n\right\rangle =\left\lfloor 1,n\right\rfloor \longrightarrow
\left\lfloor 1+1,n-1\right\rfloor =\left\lfloor 2,n-1\right\rfloor
\longrightarrow ...
\end{equation*}
\noindent ($\left\lfloor 2,n-1\right\rfloor $ is or is not equal to the
smallest set, $V$, with the property that $\left\langle n\right\rangle
\longrightarrow V$) 
\begin{equation*}
...\longrightarrow \left\lfloor 1+\left( i-2\right) ,n-\left( i-2\right)
\right\rfloor =\left\lfloor i-1,n-\left( i-2\right) \right\rfloor
\longrightarrow
\end{equation*}
\begin{equation*}
\longrightarrow \left\lfloor 1+\left( i-1\right) ,n-\left( i-1\right)
\right\rfloor =\left\lfloor i,n-\left( i-1\right) \right\rfloor
\longrightarrow
\end{equation*}
\noindent ($i$ is a left end of all the sets which follow) 
\begin{equation*}
\longrightarrow \left\lfloor i,n-i\right\rfloor \longrightarrow \left\lfloor
i,n-i-1\right\rfloor \longrightarrow ...\longrightarrow \left\lfloor
i,i+1\right\rfloor \longrightarrow \left\lfloor i,i\right\rfloor =\left\{
i\right\}
\end{equation*}
\noindent (some of the above sets can vanish --- \textit{e.g.}, $%
\left\lfloor 2,n-1\right\rfloor $ vanishes when $n=2$).

\smallskip

We must show that 
\begin{equation*}
\left\lfloor g,n-\left( g-1\right) \right\rfloor =\left\lfloor 1+\left(
g-1\right) ,n-\left( g-1\right) \right\rfloor \longrightarrow
\end{equation*}
\begin{equation*}
\longrightarrow \left\lfloor 1+g,n-g\right\rfloor ,\text{ }\forall g\in
\left\langle i-1\right\rangle ,\text{ }i\neq 1,
\end{equation*}
\noindent and 
\begin{equation*}
\left\lfloor i,h\right\rfloor \longrightarrow \left\lfloor
i,h-1\right\rfloor ,\text{ }\forall h,\text{ }i+1\leq h\leq n-\left(
i-1\right) \text{ (}i=1\text{ or }i\neq 1\text{)}
\end{equation*}
\noindent --- $i$ belongs to all the above sets.

\smallskip

Let $g\in \left\langle i-1\right\rangle ,$ $i\neq 1.$ Let $x\in \left\lfloor
g,n-\left( g-1\right) \right\rfloor .$ We show that 
\begin{equation*}
x.a\in \left\lfloor 1+g,n-g\right\rfloor
\end{equation*}
\noindent ($P=P_{a}$). Since the automaton is monotonic, we have 
\begin{equation*}
g.a\leq x.a\leq \left( n-\left( g-1\right) \right) .a.
\end{equation*}
\noindent Since $i\in \left\lfloor g,n-\left( g-1\right) \right\rfloor ,$ $%
i\neq g,$ $n-\left( g-1\right) ,$ and the automaton is monotonic, we have 
\begin{equation*}
1.a>1,\text{ }2.a>2,\text{ }...,\text{ }g.a>g,
\end{equation*}
\noindent and 
\begin{equation*}
n.a<n,\text{ }\left( n-1\right) .a<n-1,\text{ }...,\text{ }\left( n-\left(
g-1\right) \right) .a<n-\left( g-1\right)
\end{equation*}
\noindent --- indeed, 
\begin{equation*}
1.a\neq 1\Longrightarrow 1.a>1;
\end{equation*}
\begin{equation*}
2.a\neq 2\text{ and }2.a\geq 1.a>1\Longrightarrow 2.a>2;
\end{equation*}
\begin{equation*}
\vdots
\end{equation*}
\begin{equation*}
g.a\neq g\text{ (because }g\neq i\text{) and }g.a\geq \left( g-1\right)
.a>g-1\Longrightarrow g.a>g;
\end{equation*}

\medskip 
\begin{equation*}
n.a\neq n\Longrightarrow n.a<n;
\end{equation*}
\begin{equation*}
\left( n-1\right) .a\neq n-1\text{ and }\left( n-1\right) .a\leq
n.a<n\Longrightarrow \left( n-1\right) .a<n-1;
\end{equation*}
\begin{equation*}
\vdots
\end{equation*}
\begin{equation*}
\left( n-\left( g-1\right) \right) .a\neq n-\left( g-1\right) \text{
(because }n-\left( g-1\right) \neq i\text{) and }\left( n-\left( g-1\right)
\right) .a\leq
\end{equation*}
\begin{equation*}
\leq \left( n-\left( g-2\right) \right) .a<n-\left( g-2\right)
\Longrightarrow \left( n-\left( g-1\right) \right) .a<n-\left( g-1\right) .
\end{equation*}

\smallskip

Now we have 
\begin{equation*}
g<g.a\leq x.a\leq \left( n-\left( g-1\right) \right) .a<n-\left( g-1\right) .
\end{equation*}
\noindent So, 
\begin{equation*}
x.a\in \left\lfloor 1+g,n-g\right\rfloor .
\end{equation*}

\smallskip

The above proof leads, obviously, to a proof for 
\begin{equation*}
\left\lfloor i,h\right\rfloor \longrightarrow \left\lfloor
i,h-1\right\rfloor ,\text{ }\forall h,\text{ }i+1\leq h\leq n-\left(
i-1\right)
\end{equation*}
\noindent --- this proof is left to the reader.

\smallskip

By Theorem 1.7, ``$\Longleftarrow $'', 
\begin{equation*}
P^{b}
\end{equation*}
\noindent is a Markov matrix for some $b\in \left\langle n-1\right\rangle $
(we have $b=n-1$ when $i=1$).

\smallskip

\textit{Subcase} 1.2. $i-1>n-i$ ---equivalently, $i>\frac{n+1}{2}$. This
subcase is similar to Subcase 1.1 --- we have 
\begin{equation*}
\left\langle n\right\rangle =\left\lfloor 1,n\right\rfloor \longrightarrow
\left\lfloor 1+1,n-1\right\rfloor =\left\lfloor 2,n-1\right\rfloor
\longrightarrow ...
\end{equation*}
\begin{equation*}
...\longrightarrow \left\lfloor 1+n-i-1,n-\left( n-i-1\right) \right\rfloor
=\left\lfloor n-i,i+1\right\rfloor \longrightarrow
\end{equation*}
\begin{equation*}
\longrightarrow \left\lfloor 1+n-i,n-\left( n-i\right) \right\rfloor
=\left\lfloor 1+n-i,i\right\rfloor \longrightarrow
\end{equation*}
\noindent ($i$ is a right end of all the sets which follow) 
\begin{equation*}
\longrightarrow \left\lfloor 1+n-i+1,i\right\rfloor =\left\lfloor
2+n-i,i\right\rfloor \longrightarrow
\end{equation*}
\begin{equation*}
\longrightarrow \left\lfloor 3+n-i,i\right\rfloor \longrightarrow
...\longrightarrow \left\lfloor i-1,i\right\rfloor \longrightarrow
\left\lfloor i,i\right\rfloor =\left\{ i\right\} ,
\end{equation*}
\noindent so, by Theorem 1.7, ``$\Longleftarrow $'', 
\begin{equation*}
P^{b}
\end{equation*}
\noindent is a Markov matrix for some $b\in \left\langle n-1\right\rangle $
(we have $b=n-1$ when $i=n$).

\smallskip

\textit{Subcase} 1.3. $i-1=n-i$ ---equivalently, $i=\frac{n+1}{2}.$ This
subcase holds only when $n$ is odd and is also similar to Subcase 1.1 --- we
have 
\begin{equation*}
\left\langle n\right\rangle =\left\lfloor 1,n\right\rfloor \longrightarrow
\left\lfloor 1+1,n-1\right\rfloor =\left\lfloor 2,n-1\right\rfloor
\longrightarrow ...
\end{equation*}
\begin{equation*}
...\longrightarrow \left\lfloor i-1,i+1\right\rfloor \longrightarrow
\left\lfloor i,i\right\rfloor =\left\{ i\right\} ,
\end{equation*}
\noindent so, by Theorem 1.7, ``$\Longleftarrow $'', 
\begin{equation*}
P^{\frac{n-1}{2}}
\end{equation*}
\noindent is a Markov matrix.

\smallskip

From Subcases 1.1-1.3, 
\begin{equation*}
P^{b}
\end{equation*}
\noindent is a Markov matrix for some $b\in \left\langle n-1\right\rangle $
(we have $b=n-1$ when $i=1$ and when $i=n$).

\smallskip

\textit{Case} 2. $\forall P\in \mathcal{C},$ $\exists u,$ $v\in \left\langle
n\right\rangle $ ($u=u\left( P\right) ,$ $v=v\left( P\right) $)$,$ $u\neq v,$
such that $P_{uu}=P_{vv}=1.$ In this case, we do not find any proof using
the $G^{+}$ method yet. Note that, in this case, any Markov product, if any,
contains at least two different matrices from $\mathcal{C}.$ We cannot have
Markov products (of matrices from $\mathcal{C}$) if $\left| \mathcal{C}%
\right| =1$ or, more generally, $\exists i,$ $j\in \left\langle
n\right\rangle ,$ $i\neq j,$ such that $P_{ii}=P_{jj}=1,$ $\forall P\in 
\mathcal{C}$. An example in which we have Markov products: $\mathcal{C}%
=\left\{ P_{1},P_{2}\right\} ,$%
\begin{equation*}
P_{1}=\left( 
\begin{array}{ccc}
1 & 0 & 0 \\ 
0 & 1 & 0 \\ 
0 & 1 & 0%
\end{array}
\right) ,\text{ }P_{2}=\left( 
\begin{array}{ccc}
1 & 0 & 0 \\ 
1 & 0 & 0 \\ 
0 & 0 & 1%
\end{array}
\right)
\end{equation*}
\noindent --- we have 
\begin{equation*}
P_{1}P_{2}=\left( 
\begin{array}{ccc}
1 & 0 & 0 \\ 
1 & 0 & 0 \\ 
1 & 0 & 0%
\end{array}
\right) ,
\end{equation*}
\noindent so $P_{1}P_{2}$ is a Markov matrix/product ($\Longrightarrow
P_{1}P_{2}P_{1},$ $P_{1}P_{2}P_{2},$ $P_{1}P_{2}P_{1}P_{1},$ etc. are Markov
matrices), and the $G^{+}$ method can be applied for $P_{1}P_{2}$.

\smallskip

The $G^{+}$ method can be used, on the \v{C}ern\'{y} conjecture, both to
find the upper bounds (see Examples 2.19 and 2.21 and Remark 2.20) and to
find reset words having as small as possible length (see Example 2.19...).

\smallskip

For more information on the synchronizable deterministic finite automata,
for applications and the \v{C}ern\'{y} conjecture, see, \textit{e.g.}, [5],
[9]-[10], [12]-[16], [19], [30] ([15] is in connection with [30]), [34],
[37]-[40], and [42] --- for applications, see [39]...

\bigskip

\begin{center}
\textbf{REFERENCES}
\end{center}

\bigskip

[1] D.S. Ananichev and M.V. Volkov, \textit{Synchronizing monotonic automata}%
. Theoret. Comput. Sci. \textbf{327} (2004), 225--239.

[2] J. Bang-Jensen and G.Z. Gutin, \textit{Digraphs}:\textit{\ Theory,
Algorithms and Applications}, 2nd Edition. Springer, London, 2009.

[3] R. Bellman, \textit{Introduction to Matrix Analysis}. McGraw-Hill Book
Company, New York, 1960.

[4] A. Benjamin, G. Chartrand, and P. Zhang, \textit{The Fascinating World
of Graph Theory}. Princeton University Press, Princeton, 2015.

[5] M.V. Berlinkov and M. Szyku\l a, \textit{Algebraic synchronization
criterion and computing reset words}. Inform. Sci. \textbf{369} (2016),
718--730.

[6] V.D. Blondel and A. Olshevsky, \textit{How to decide consensus}? \textit{%
A combinatorial necessary and sufficient condition and a proof that
consensus is decidable but NP-hard}. SIAM J. Control Optim. \textbf{52}
(2014), 2707-2726.

[7] V.D. Blondel, R.M. Jungers, and A. Olshevsky, \textit{On primitivity of
sets of matrices}. Automatica J. IFAC \textbf{61} (2015), 80--88.

[8] J.A. Bondy and U.S.R. Murty, \textit{Graph Theory with Applications}.
North-Holland, New York, 1979.

[9] J. \v{C}ern\'{y}, \textit{Pozn\'{a}mka k homog\'{e}nnym experimetom s
kone\v{c}n\'{y}mi automatmi}. Matematicko-fyzik\'{a}lny \v{C}asopis, Slovensk%
\'{a} Akad\'{e}mia Vied \textbf{14} (1964), 208--216. (In Slovak; English
translation: \textit{A note on homogeneous experiments with finite automata}%
, J. Autom. Lang. Comb. \textbf{24} (2019), 123--132.)

[10] P.-Y. Chevalier, J.M. Hendrickx, and R.M. Jungers, \textit{Reachability
of consensus and synchronizing automata}. 2015 54th IEEE Conference on
Decision and Control (CDC), 4139--4144, 2015.

[11] P.-Y. Chevalier, V.V. Gusev, J.M. Hendrickx, and R.M. Jungers, \textit{%
Sets of stochastic matrices with converging products}:\textit{\ bounds and
complexity}. arXiv:1712.02614v1.

[12] H. Don, \textit{The \v{C}ern\'{y} conjecture and }1\textit{-contracting
automata}. Electron. J. Combin. \textbf{23(3)} (2016), \#P3.12 (10 pages).

[13] D. Eppstein, \textit{Reset sequences for monotonic automata}. SIAM J.
Comput. \textbf{19} (1990), 500--510.

[14] D. Eppstein and J.-C. Falmagne, \textit{Algorithms for media}. Discrete
Appl. Math. \textbf{156} (2008), 1308--1320.

[15] P. Frankl, \textit{An extremal problem for two families of sets}.
European J. Combin. \textbf{3} (1982), 125--127.

[16] B. Gerencs\'{e}r, V.V. Gusev, and R.M. Jungers, \textit{Primitive sets
of nonnegative matrices and synchronizing automata}. SIAM J. Matrix Anal.
Appl. \textbf{39} (2018), 83--98.

[17] D.J. Hartfiel, \textit{Nonhomogeneous Matrix Products}. World
Scientific, Singapore, 2002.

[18] M. Iosifescu, \textit{Finite Markov Processes and Their Applications}.
Wiley, Chichester \& Ed. Tehnic\u{a}, Bucharest, 1980; corrected
republication by Dover, Mineola, N.Y., 2007.

[19] H. J\"{u}rgensen, \textit{Synchronization}. Inform. and Comput. \textbf{%
206} (2008), 1033--1044.

[20] D.E. Knuth, \textit{The Art of Computer Programming}: Vol. \textbf{4A}, 
\textit{Combinatorial Algorithms}, Part 1. Addison-Wesley, Boston, 2011.

[21] H. Minc, \textit{Permanents}. Cambridge University Press, Cambridge,
1984.

[22] U. P\u{a}un, $\Delta $\textit{-ergodic theory and reliability theory}.
Math. Rep. (Bucur.) \textbf{10}(\textbf{60}) (2008), 73$-$95.

[23] U. P\u{a}un, $G_{\Delta _{1},\Delta _{2}}$ \textit{in action}. Rev.
Roumaine Math. Pures Appl.\textbf{\ 55 }(2010), 387$-$406.

[24] U. P\u{a}un, \textit{A hybrid Metropolis-Hastings chain}. Rev. Roumaine
Math. Pures Appl.\textbf{\ 56 }(2011), 207$-$228.

[25] U. P\u{a}un, $P\left( X_{s}\in A_{s},X_{s+1}\in A_{s+1},...,X_{t}\in
A_{t}\right) $ \textit{in the Markov chain case}:\textit{\ from an upper
bound to a method}. Rev. Roumaine Math. Pures Appl.\textbf{\ 57 }(2012), 145$%
-$158.

[26] U. P\u{a}un, \textit{Waiting time random variables}: \textit{upper
bounds}. Markov Process. Related Fields \textbf{19} (2013), 791$-$818.

[27] U. P\u{a}un, \textit{Other results on the Markovian inequality}

\noindent $P\left( X_{s}\in A_{s},X_{s+1}\in A_{s+1},...,X_{t}\in
A_{t}\right) \leq \overset{\_}{\alpha }\left( Q_{s,t}\right) .$ Rev.
Roumaine Math. Pu-

\noindent res Appl.\textbf{\ 61 }(2016), 157$-$183.

[28] U. P\u{a}un, $G^{+}$ \textit{method in action}:\textit{\ new classes of
nonnegative matrices}, \textit{with results.} ROMAI\ J. \textbf{20} (2024),
no. 1, 99$-$143. With preprint: arXiv:2304.05227.

[29] U. P\u{a}un, $G$\textit{\ method and finite-time consensus.} ROMAI\ J. 
\textbf{22} (2026), no. 1, 89$-$127. With preprint: arXiv:2510.17542.

[30] J.-\'{E}. Pin, \textit{On two combinatorial problems arising from
automata theory}. North-Holland Math. Stud. \textbf{75} (1983), 535--548 or,
but only in printed format, Ann. Discrete Math. \textbf{17} (1983), 535--548.

[31] V. Yu. Protasov and A.S. Voynov, \textit{Sets of nonnegative matrices
without positive products}. Linear Algebra Appl. \textbf{437} (2012),
749--765.

[32] V. Yu. Protasov and R.M. Jungers, \textit{Lower and upper bounds for
the largest Lyapunov exponent of matrices}. Linear Algebra Appl. \textbf{438}
(2013), 4448--4468.

[33] M.O. Rabin, \textit{Probabilistic Automata}. Information and Control 
\textbf{6} (1963), 230--245.

[34] I. Rystsov, \textit{Reset words for commutative and solvable automata}.
Theoret. Comput. Sci. \textbf{172} (1997), 273--279.

[35] E. Seneta, \textit{Non-negative Matrices and Markov Chains}, 2nd
Edition. Springer-Verlag, Berlin, 1981; revised printing, 2006.

[36] Y. Shitov, \textit{An improvement to a recent upper bound for
synchronizing words of finite automata}. J. Autom. Lang. Comb. \textbf{24}
(2019), 367--373.

[37] A. N. Trahtman, \textit{The \v{C}ern\'{y} conjecture for aperiodic
automata}. Discrete Math. Theoret. Comput. Sci. \textbf{9} (2007), 3--10.

[38] M.V. Volkov, \textit{Synchronizing automata and the \v{C}ern\'{y}
conjecture.} In: C. Mart\'{\i}n-Vide, F. Otto, and H. Fernau (Eds.), LATA
2008, LNCS \textbf{5196} (2008), 11--27, Springer-Verlag, Berlin.

[39] M.V. Volkov, \textit{Synchronization of finite automata}. Russian Math.
Surveys \textbf{77} (2022), 819--891.

[40] M.V. Volkov, \textit{List of results on the \v{C}ern\'{y} conjecture
and reset thresholds for synchronizing automata.} arXiv:2508.15655v3.

[41] A. Voynov, \textit{Shortest positive products of nonnegative matrices}.
Linear Algebra Appl. \textbf{439} (2013), 1627--1634.

[42] Y. Zhu, \textit{A quadratic upper bound on the reset thresholds of
synchronizing automata containing a transitive permutation group}. In: S.
Barman and S. Lasota (Eds.), 44th IARCS Annual Conference on Foundations of
Software Technology and Theoretical Computer Science (FSTTCS 2024), LIPICS,
Vol. \textbf{323}, pp. 34:1--34:14, Schloss Dagstuhl - Leibniz-Zentrum f\"{u}%
r Informatik, 2024.

\bigskip\ 
\begin{equation*}
\begin{array}{ccc}
\begin{array}{c}
\\ 
\mathit{December}\text{ }\mathit{1,}\text{ }\mathit{2026}%
\end{array}
&  & \mathit{Romanian\ Academy} \\ 
&  & \mathit{Gheorghe\ Mihoc}-\mathit{Caius\ Iacob}\text{ }\mathit{Institute}
\\ 
&  & \mathit{of\ Mathematical\ Statistics}\text{ }\mathit{and\ Applied\
Mathematics} \\ 
&  & \mathit{Calea\ 13\ Septembrie\ nr.\ 13} \\ 
&  & \mathit{050711\ Bucharest\ 5,\ Romania} \\ 
&  & \mathit{upterra@gmail.com}%
\end{array}%
\end{equation*}

\end{document}